\UseRawInputEncoding
\documentclass[12pt]{amsart}
\usepackage{amsmath}\usepackage{MnSymbol}
\usepackage{amsfonts}\usepackage{verbatim}\usepackage{bbold}
\usepackage{graphicx}\usepackage{fontawesome}
\usepackage{marvosym}\usepackage{ifsym}
\usepackage{dictsym}
\usepackage{wasysym}
\usepackage{amstext}
\usepackage{amsbsy}
\usepackage{amsopn}\usepackage{hyperref}
\usepackage{amsthm}\usepackage{color}

\theoremstyle{definition}

\theoremstyle{remark}

\def\proclaim#1{\vskip0.5em\noindent{\bf #1}\it}
\def\endproclaim{\vskip0.5em\par\noindent\rm}
\def\proclaim#1{\vskip0.5em\noindent{\bf #1}\it}
\def\endproclaim{\vskip0.5em\par\noindent\rm}
\def\demo#1{\vskip0.5em\noindent{\bf #1\ }}

\def\text#1{\mbox{#1}}
\def\flushpar{\par\noindent}

\def\tag#1{\eqno{(#1)}}
\def\mod{\mbox{ mod }}
\newcommand{\mapright}[1]{%
    \smash{\mathop{%
        \hbox to 1cm{\rightarrowfill}
        }
    \limits^{#1}
    }
}
\newcommand{\mapleft}[1]{%
    \smash{\mathop{%
        \hbox to 1cm{\rightarrowfill}
        }
    \limits_{#1}
    }
}

\usepackage{wasysym}\usepackage{eiad}
\def\e{\epsilon}
\def\a{\alpha}
\def\b{\beta}
\def\g{\gamma}
\def\G{\Gamma}
\def\d{\delta}
\def\D{\Delta}
\def\th{\theta}
\def\s{\sigma}

\def\l{\lambda}
\def\L{\Lambda}
\def\x{\times}

\def\o{\overline}
\def\f{\flushpar}

\def\v{\varphi}

\def\bdy{\partial}

\def\om{\omega}
\def\Om{\Omega}
\def\B{\Cal B}

\def\A{\Cal A}

\def\T{\widehat T}
\def\({\biggl(}
\def\){\biggr)}
\def\[{\bigl[}
\def\]{\bigr]}\def\Cal{\mathcal}

\def\<{\langle}
\def\>{\rangle}

\def\bul{\smallskip\f$\bullet\ \ \ $}\def\Par{\smallskip\f\P}
\def\sms{\smallskip\f}\def\sbul{\f$\bullet\
\ \ $}\def\sms{\smallskip\f}
\def\st{\text{such that}}\def\lra{\longrightarrow}

\def\Smi{\smallskip\f{\LARGE \smiley\ \ \ }}

\def\xbm{(X,m)}
\def\lcl{\lceil}\def\rcl{\rceil}\def\ttau{\widehat{\tau}}
\begin{document}

 \title{
 tied-down  occupation times of infinite ergodic transformations}
\author{ Jon. Aaronson and Toru Sera}
\address[Aaronson]{School of Math. Sciences, Tel Aviv University,
69978 Tel Aviv, Israel.}
\email{aaro@tau.ac.il}
\address[Sera]{Department of Mathematics, Graduate School of Science, Osaka University,
Toyonaka, Osaka 560-0043, Japan}
\email{ sera@cr.math.sci.osaka-u.ac.jp}

\begin{abstract}We prove distributional limit theorems (conditional and integrated) for the occupation times  of certain  weakly mixing,
pointwise dual ergodic transformations at  ``tied-down'' times immediately after ``excursions''.
The limiting random variables include the local times of $q$-stable L\'evy-bridges ($1<q\le 2$) and
    the transformations involved exhibit  ''tied-down renewal mixing`` properties which refine   rational weak mixing. 
    Periodic local limit theorems for Gibbs-Markov and AFU maps are also established.
\end{abstract}
\thanks{\copyright 2019-20. Aaronson{\tiny$^\prime$}s research was partially supported by ISF grant No. 1289/17.
Sera{\tiny$^\prime$}s research was partially supported by JSPS KAKENHI grants No. 19J11798 and No. JP21J00015.}
\maketitle\markboth{tied-down  occupation times}{Jon. Aaronson and Toru Sera}

\section*{\S0 Introduction}
Let $(X,m,T)$   denote   a measure preserving transformation ({\tt MPT}) $T$ of the non-atomic, Polish measure space $(X,m)$ (where 
$m$ is a $\s$-finite, non-atomic measure
defined  on the Borel subsets $\B(X)$ of $X$). 
\

The associated {\it transfer operator} $\T:L^1(m)\to L^1(m)$ is the predual of $f\mapsto f\circ T\ \ (f\in L^\infty(m))$, that is 
$$\int_X\T f gdm=\int_Xf g\circ Tdm\ \ \text{for}\ f\in L^1(m)\ \&\ g\in L^\infty(m).$$
\

The {\tt MPT} $(X,m,T)$ is called  {\it pointwise dual ergodic} ({\tt PDE}) if there is a sequence $a(n)=a_n(T)$ (the {\it return sequence} of $(X,m,T)$) so that
\begin{align*}\tag{{\tt PDE}}\label{PDE}\frac1{a(n)}\sum_{k=0}^{n-1}\T^kf\xrightarrow[n\to\infty]{}\ \int_Xfdm\ \text{a.e.}\ \forall\ f\in L^1(m). 
\end{align*}

Pointwise dual ergodicity entails 
\sbul {\tt\small recurrence} (conservativity-- no non-trivial wandering sets), and 

\sbul{\tt\small ergodicity} (no non-trivial invariant sets).
\

\subsection*{Distributional limits}
\

Let  $(X,m,T)$ be   pointwise dual ergodic with $\g$-regularly varying return sequence $a(n)=a_n(T)$ ($0<\g<1$). 
\

By the Darling Kac theorem (\cite{DK})\footnote{see also the  version in  \cite[Ch. 3]{IET}}
\begin{align*}\tag{{\scriptsize\faAnchor}}\label{faAnchor}
 \tfrac1{a(n)}S_n(f)\overset{\mathfrak d}\lra \ Y_\g m(f)\ \forall\ f\in L^1_+
\end{align*}

where
$Y_\g$ is the Mittag-Leffler distribution of order $\g$ (see \cite[XIII.8]{Feller}) normalized so that $\Bbb E(Y_\g)=1$,  $m(f):=\int_Xfdm$ and $\overset{\mathfrak d}\lra$ (on $(X,m)$)
denotes convergence in distribution with respect to  $m$-absolutely continuous probabilities.

Let $\Om\in\mathcal F_+:=\{A\in B(X):\ \ 0<m(A)<\infty\}$. The {\it return time function} to $\Om$ is $\v=\v_\Om:\Om\to\mathbb N$ defined by
$\v(\om):=\min\{n\ge 1:\ T^n\om\in\Om\}<\infty$ a.s. by conservativity. 
\

\subsection*{Return time process}
\

The {\it induced transformation} on  $\Om$ is $T_\Om:\Om\to\Om$ defined by $T_\Om(\om):=T^{\v(\om)}(\om)$. As shown in
\cite{Kakutani}, it is an ergodic, probability preserving transformation of $(\Om,m_\Om)$ (where $m_\Om(A):=\tfrac{m(\Om\cap A)}{m(\Om)}$).
\

A standard inversion argument on (\ref{faAnchor}) 
shows that the {\it return time process}\label{retproc} 
$(\Om,m_\Om,T_\Om,\v_\Om)$ on $\Om$ satisfies the {\emph{ stable limit theorem}}:
\begin{align*}\tag{{\tt SLT}}\label{SLT}
  \frac{\v_n}{a^{-1}(n)}\xrightarrow[n\to\infty]{\mathfrak d}\ \  \frac{Z_\g}{m(\Om)^{1/\g}}
\end{align*}
  on $(\Om,m_\Om)$ where
 $\v_n:=\sum_{k=0}^{n-1}\v\circ T_\Om^k$.
 \
 
 As shown in \cite[XIII.8]{Feller}, $Z_\g:=Y_\g^{-\frac1\g}$ is the positive, stable random variable of order $\g$ with Laplace transform 
 $\Bbb E(e^{-sZ_\g})=\exp[-\tfrac{s^\g}{\G(1+\g)}]$\ ($s>0$) whence with
 characteristic function 
 $$\Phi_{Z_\g}(t):=\Bbb E(e^{it Z_\g})=\exp[-\tfrac{|t|^\g}{\G(1+\g)}(\cos\tfrac{\g\pi}2-i\text{sgn}(t)\sin\tfrac{\g\pi}2)].$$
\

The stable limit theorem holds for the return time process of any $\Om\in\mathcal F_+$. 
By ``choosing" $\Om$ carefully, it is sometimes possible to obtain   stronger properties for the return time process and  possibly also for the {\tt MPT}.
\

\subsection*{Tied-down, renewal mixing properties}
 \
In this paper we study additional properties of  certain  pointwise dual ergodic  transformations
  preserving infinite  measures. 
  Analogous properties for invertible {\tt MPT}s  are considered
  in \S6. 
  
  \
  
  Our additional properties are related to ''{\tt\small tied-down renewal theory}`` which studies renewal counts at renewal times as   in \cite{Wendel}, \cite{Liggett1970} and \cite{Godreche} and extend the strong renewal theorem
  as in \cite{G-L, Melter,Goue11, C-D}.
  \

 We'll consider the following  conditional, {\tt\small tied-down renewal mixing properties} of  a {\tt PDE MPT}   $(X,m,T)$
with $a(n)=a_n(T)$  
$\g$-regularly varying   ($\g\in (0,1)$):
\begin{align*}\tag{\dsmilitary}\label{dsmilitary}\tfrac1{a(N)}\sum_{n=1}^N&|\T^n(1_Cg(\tfrac{S_n(f)}{a(n)}))-m(C) \Bbb E(g(m(f)W_\g))u(n)|
\xrightarrow[N\to\infty]{}\ 0\\ & \text{a.e.}\ \forall\ C\in\mathcal F_+,\ g\in C_B(\Bbb R_+)\ \&\ f\in L^1_+ 
\end{align*}
where $u(n)\sim\tfrac{\g a(n)}N,\ W_\g\in\text{\tt RV}\,(\Bbb R_+),\  \Bbb E(g(W_\g))=\Bbb E(Y_\g g(Y_\g))$.
\

The convergence mode in (\dsmilitary) is considered in \cite[\S4]{RatWM} as $u$-{\tt\small  strong Cesaro convergence}.
Direct convergence is impossible e.g. for $C\in\mathcal{F}_+$ a {\tt\small weakly wandering set} as in \cite{Haj-Kak1}.

 \
 
 The limit random variables $W_\g$ ($0<\g<1$) appear in \cite{Wendel}. For $0<\g\le\frac12$ they correspond to the local time at zero up to time $1$ of the {\tt\small symmetric $\frac1{1-\g}$-stable bridge} and for $0<\g<1$ to the local time at zero up to time $1$ of the {\tt\small  Bessel bridge of dimension} $2-2\gamma$ (see \cite{PitYor97}). 
\

Note that (\ref{dsmilitary}) \ entails conditional rational weak mixing as in \cite{A-T}, whence spectral weak mixing and ergodicity of $(X,m,T^N)\ \forall\ N\ge 1$. 
Calculation now shows that if $(X,m,T)$ satisfies (\dsmilitary), then so does $(X,m,T^N)\ \forall\ N\ge 1$ with $a_n(T^N)\sim\tfrac{a_{nN}(T)}N$.

\subsection*{Local limits for asymptotically  stable, stationary processes }\label{assp}
For $\g\in (0,1)$, call the positive, ergodic, stationary process 
$(\Om,\mu,\tau,\phi)$ (i.e. where $(\Om,\mu,\tau)$ is an ergodic probability preserving transformation ({\tt EPPT}) $\&\ \phi:\Om\to\Bbb R$ is measurable) {\it asymptotically $\g$-stable} \  if 
\begin{align*}\tag{{\scriptsize\faKey}}\label{faKey}
 \mu([\phi>t])\  \sim\ \frac1{\G(1+\g)\G(1-\g)}\cdot\frac1{a(t)}.
\end{align*}
 with $a(t)\ \ \ \g$-regularly varying and $(\Om,\mu,\tau,\phi)$  satisfies  the stable limit theorem as in (\ref{SLT}) on p. \pageref{SLT}.
 \

Let $(\Om,\mu,\tau,\phi)$ be an  asymptotically $\g$-stable, positive, ergodic, stationary process with
$\phi:X\to\Bbb G=\Bbb Z$ in the lattice case or $\Bbb G=\Bbb R$ in the continuous case.
\

Let $\a\subset\B(\Om)$ be a countable partition  so that $(\Om,\mu,\tau,\a)$ is a {\tt\small fibered system} as defined on p. \pageref{fibered}.
\

Denote $\a_n:=\bigvee_{k=0}^{n-1}\tau^{-k}\a\ \&\ \mathcal C_\a:=\bigcup_{n\ge 1}\a_n$,
the collection of $\a$-{\it cylinders}.

\

We'll say that $(\Om,\mu,\tau,\a,\phi)$ satisfies 
\bul an {\it aperiodic local limit theorem}\label{apLLT} ({\tt LLT}) if $\forall\ \ A\in\mathcal C_\a,\ I\subset \Bbb G$ a bounded interval,
\begin{align*}\tag{\Coffeecup}\label{Coffeecup} a^{-1}(n)\widehat{\tau}^n (1_{A\cap  [\phi_n\in k_n+I]})
\xrightarrow[n\to\infty,\ k_n\in\Bbb G,\ \tfrac{k_n}{a^{-1}(n)}\to x]{}\ f_{Z_\g}(x)m_\Bbb G(I)m(A)\end{align*}
uniformly in $x\in [c,d]$ whenever $0<c<d<\infty$, where $m_\Bbb G$ denotes Haar measure on $\Bbb G$ and$f_{Z_\g}$ is the probability density function of $Z_\g$;

\bul a {\it periodic {\tt LLT}} \label{LLT}if $\exists\ p\in\Bbb G,\ p>0$ (the {\it period})\ 
and $\xi\in (0,p)\cap\Bbb G$ (the {\it drift}) so that

\

$\overline{\<\{n\xi:\ n\ge 1\}+p\Bbb Z\>}=\Bbb G$ 
  $\&\ \forall\ \ A\in\mathcal C_\a,\ I\subset \Bbb G$ a bounded interval, with $I_J:=I\cap [Jp,(J+1)p)$ ($J\in\Bbb Z$),

\begin{align*}\tag{\dsrailways}\label{dsrailways} a^{-1}(n)\widehat{\tau}^n (1_{A\cap  [\phi_n\in k_n+I]})&
\approx\  pf_{Z_\g}(x)m(A)\sum_{J\in\Bbb Z}1_{p\Bbb Z+I_J}(k_n-n\xi)\\ &\text{as}\ \ n\to\infty,\ k_n\in\Bbb G,\ \tfrac{k_n}{a^{-1}(n)}\to x \end{align*}
uniformly in $x\in [c,d]$ whenever $0<c<d<\infty$;
\

and
\sbul\ a {\it generalized {\tt LLT}} ({\tt GLLT}) if it satisfies either an aperiodic, or a periodic {\tt LLT}.
\

Here and throughout, $a_n\approx b_n$ means $a_n-b_n\xrightarrow[n\to\infty]{}0$.
\

Local limits are studied (here) using Fourier theory and throughout the paper,\label{char}
 for  $\Bbb H$ a locally compact, Polish, Abelian group, $\widehat{\Bbb H}$ denotes the multipicative group of {\it characters}, that is
$$\widehat{\Bbb H}:=\{\chi:\Bbb H\to\Bbb S^1:\ \chi\ \text{a continuous homomorphism}\}$$
where $\Bbb S^1:=\{z\in\Bbb C:\ |z|=1\}$. 
\

For example $\widehat{\Bbb R^d}\cong\Bbb R^d$ via $\chi_t(x)=e(\<x,t\>)$ ($x,t\in\Bbb R^d$) where $\<x,t\>:=\sum_{j=1}^dx_jt_j$ and
$e(s):=e^{2\pi is}$ ($s\in\Bbb R$).
\subsection*{Local limits for periodic cocycles}\label{Periodicity}
\

\

Let $(\Om,\mu,\tau)$ be an ergodic, probability preserving transformation and let $\Bbb G=\Bbb Z$ or $\Bbb R$.  
\

The measurable function $\phi:\Om\to \Bbb G$ is called:
\bul {\it non-arithmetic} if for  $\chi\in\widehat{\Bbb G},\ \chi\nequiv 1$, $\chi\circ\phi$ is not a $\tau$-coboundary, i.e. the equation 
$\chi\circ\phi=\tfrac{g(x)}{g(Tx)}$ where $g:\Om\to\Bbb S^1$ is measurable, has no  solutions; and 
\bul  {\it aperiodic } if for  $\chi\in\widehat{\Bbb G},\ \chi\nequiv 1,\ \l\in\Bbb S^1$,
$\l\chi\circ\phi$ is not a $\tau$-coboundary.

\

It is not hard to check that if $(\Om,\mu,\tau)$ is weakly mixing, and $\phi:\Om\to\Bbb G$, then 
    \sms (i) in case the skew product $(\Om\x\Bbb G,\mu\x m_\Bbb G,\tau_\phi)$ is ergodic, it is  weakly mixing iff $\phi$ is aperiodic; and
\sms (ii) in case $\phi>0$, $(\Om,\mu,\tau)^\phi$
  is weakly mixing iff $\phi$ is non-arithmetic. 
  \
  
 Here, in case $\Bbb G=\Bbb Z$, $(\Om,\mu,\tau)^\phi$ denotes  the {\it Kakutani skyscraper}\label{skyscraper}:
\

the conservative, ergodic {\tt MPT} ({\tt CEMPT} ) defined in  \cite{Kakutani} by  $(\Om,\mu,\tau)^\phi:=(X,m,T)$ where
\begin{align*}
 X:=\{(\om,n)\in\Om\x\Bbb N:\ &1\le n\le \phi(\om)\},\ m:=\mu\x\#|_X\ \&\\ & T(\om,n):=\begin{cases}& (\om,n+1)\ \ n<\phi(\om)\\ &
                                                                                   (\tau(\om),1)\ \ n=\phi(\om)
                                                                                  \end{cases}
\end{align*}
 with return time process to $\Om\x\{1\}$ 
$\cong (\Om,\mu,\tau,\phi).$
 \
 
 In case $\Bbb G=\Bbb R$, $(\Om,\mu,\tau)^\phi$ denotes   the  suspended semiflow  over $(\Om,\mu,\tau)$ under the {\it ceiling } $\phi$ (p. \pageref{semiflow}). 
 \

We'll establish in \S3 the  periodic {\tt LLT} (Theorem 3.1 on p. \pageref{GenLLT})\   for  $(\Om,\mu,\tau,\a)$ a
{\tt D-F} fibered system (p. \pageref{DFfibered}) which is Gibbs-Markov or {\tt AFU} and 
$\phi:\Om\to\Bbb G$ a non-arithmetic, {\tt D-F} cocycle (p. \pageref{DFcocycle}). This uses the aperiodic {\tt LLT}s of \cite{AD,ADSZ}.  
\ 

Periodic {\tt LLT}s for independent, identically distributed random variables were established in \cite{Shepp}.

\subsection*{Local limit sets}
\

\ Let $(X,m,T)$ be a {\tt CEMPT}.
\

We'll say that $\Om\in\B(X),\ 0<m(\Om)<\infty$ is a 
{\it  generalized local limit}  ({\tt GLL})  {\it set for $T$}  if
the return time process $(\Om,m_\Om,T_\Om,\v_\Om)$ satisfies a  lattice {\tt GLLT} 
with respect to some one-sided, countable generator $\a\subset\B(\Om)$. In \S5, we'll consider {\tt\small semiflows} with {\tt\small sections} 
satisfying continuous {\tt GLLT}s.

 \section*{\S1 Results $\&$ examples}
 \
 
 \subsection*{Tied down renewal mixing}
  \
  
  \ 
 Our main result is
 \
 
 \proclaim{Theorem A}
\

Suppose that  $(X,m,T)$ is pointwise dual ergodic with $a(n)=a_n(T)$  
$\g$-regularly varying {\rm($\g\in (0,1)$)} and which, has a  {\tt GLL}  set $\Om\in\mathcal F_+$, then  $(X,m,T)$ satisfies (\ref{dsmilitary})  {\rm (on p. \pageref{dsmilitary})}.
\endproclaim
 In particular  
 $(X,m,T)$ is  conditionally rationally weakly mixing. This latter property also follows, in the aperiodic case, from proposition 3.1 of \cite{nice}.
 
 \

 We'll prove theorem A in  \S2. 
 In \S5 we consider analogous properties for semiflows.
  Theorem 5.2  is a continuous time version of theorem A. 
   \subsection*{Existence of {\tt GLL} sets}
 \
 
 Asymptotically $\g$-stable, positive, ergodic, stationary, stochastic processes  are considered in \S3,
 where sufficient conditions for their {\tt GLLT}s are established via:
 
  \proclaim{Theorem B}
 \
 
 \ Suppose that $(X,m,T)$ is a pointwise dual ergodic, weakly mixing {\tt MPT} with $a(n)=a_n(T)\ \g$-regularly 
 varying with $0<\g<1$.
 \
 
 Suppose that the return time process on $\Om\in\B(X),\ 0<m(\Om)<\infty$ admits a one-sided generator $\a\subset\B(\Om)$
 so that  $(\Om,m_\Om,T_\Om,\a)$ is a mixing and either a Gibbs-Markov map or an {\tt AFU} map and $\v_\Om$ is $\a$-measurable, then 
 $\Om$ is a {\tt GLL} set for $T$.\endproclaim
 \subsection*{Examples:\ Tied-down renewals and the strong renewal theorem}
\ \ 

\

Let $(\Om,\mu,\tau,\v)$ be a positive, $\Bbb G$-valued,  Bernoulli process\footnote{i.e. where $(\v\circ\tau^n:\ n\ge 0)$ are iidrvs} with $\Bbb G=\Bbb Z$ ({\it discrete time} case) or $\Bbb R$ ({\it continuous time} case)
which is asymptotically $\g$-stable  with $\g\in (0,1)$ (i.e. satisfies (\ref{faKey}) as on p. \pageref{faKey}). 
\

Suppose that  the support of $\v$ generates $\Bbb G$ in the sense that
$$\overline{\<\{y\in\Bbb G:\ \mu([|\v-y|<\e])>0\ \forall\ \e>0\}\>}=\Bbb G.$$
\ \ 

The following results introduce the ideas of theorems A and 5.2  respectively.
\subsubsection*{\tt Tied down  renewals: discrete time}\ \ 
\

 In case $\Bbb G=\Bbb Z$,
 for $g\in C_B(\Bbb R_+)$ and with 
$u(n):=\frac{\g a(n)}n$, there is a set $K\subset\Bbb N$ with zero density ($|K\cap [1,n]|=o(n)$) so that
\begin{align*}\tag{{\scriptsize\faLightbulbO}}\label{faLightbulbO}
\frac1{u(n)}\sum_{k=1}^ng(\tfrac{k}{a(n)})\mu([\v_k=n])\xrightarrow[n\to\infty,\ n\notin K]{×} \Bbb E(g(W_\g))
\end{align*}
where $W_\g$ is as in Theorem A.

\subsubsection*{\tt Tied down renewals: continuous time}\ \
\

 In case $\Bbb G=\Bbb R$, then  for $g\in C_B(\Bbb R_+)$ and $I\subset (0,1)$ an interval; with 
$u(n):=\frac{\g a(n)}n$,  there is a set $K\subset\Bbb N$ with zero density ($|K\cap [1,n]|=o(n)$) so that
\begin{align*}\tag{{\scriptsize\faBatteryQuarter}}\label{faBatteryQuarter}\frac1{u(n)}\sum_{k\ge 1}g(\tfrac{k}{a(n)})\mu([\v_k\in n+I])\xrightarrow[n\to\infty,\ n\notin K]{×}|I| \Bbb E(g(W_\g))\big|.
\end{align*}

\

\demo{Sketch proof of (\ref{faLightbulbO})}
\

By \cite{Shepp}, $\Om$ is a {\tt GLL} set 
for   $(X,m,T)=(\Om,\mu,\tau)^\v$, the Kakutani skyscraper as on p. \pageref{skyscraper}.  By Lemma 2.1 below (on p. \pageref{lemma2.1}),
\begin{align*}\frac1{a(N)}\sum_{n=1}^N
\big|\sum_{k=1}^ng(\tfrac{k}{a(n)})\mu([\v_k=n])- \Bbb E(g(W_\g))u(n)\big|\xrightarrow[N\to\infty]{}\ 0
\end{align*}

whence  (see e.g. \cite{G-L}) (\ref{faLightbulbO}).\ \Checkedbox

\

\demo{Sketch proof of (\ref{faBatteryQuarter})}
\

Also by \cite{Shepp}, $(\Om,\mu,\tau,\v)$  is a {\tt GLL} section 
for  the suspension semiflow $(X,m,\Psi)=(\Om,\mu,\tau)^\v$ (see \S5). In case it is not a {\tt\small standard section} (as on p.\pageref{standard}), apply
 remark 5.5 (on p. \pageref{remark5.5}) to obtain a suitable  Bernoulli induced section, which, by \cite{Shepp}, is a {\tt GLL} section. Theorem  5.2  (on p. \pageref{theorem5.2}) with $t=1$ applies and 
\begin{align*}\frac1{a(N)}\sum_{n=1}^N
\big|\sum_{k\ge 1}g(\tfrac{k}{a(n)})\mu([\v_k\in n+I])-u(n)|I| \Bbb E(g(W_\g))\big|\xrightarrow[N\to\infty]{}\ 0.
\end{align*}
whence  (see e.g. \cite{G-L}) (\ref{faBatteryQuarter}).\ \Checkedbox
\

Convergence in (\ref{faLightbulbO}) and (\ref{faBatteryQuarter}) is equivalent to the {\tt\small strong renewal theorem} (SRT) as in \cite{C-D}.
\

In the discrete case,
\begin{align*}&\frac1{u(n)}\sum_{k=1}^ng(\tfrac{k}{a(n)})\mu([\v_k=n])\xrightarrow[n\to\infty]{}\  \Bbb E(g(W_\g))\ \ \forall\ g\in C_B(\Bbb R_+)\\ &\  
\iff\ \\ &\tag{\faBolt}\sum_{k=1}^n\mu([\v_k=n])\sim u(n);
\end{align*}
and in the continuous case,
\begin{align*}&\frac1{u(n)}\sum_{k=1}^ng(\tfrac{k}{a(n)})\mu([\v_k=n+I])\xrightarrow[n\to\infty]{}\ |I| \Bbb E(g(W_\g))\\ & \forall\ g\in C_B(\Bbb R_+)\ \&\ \text{intervals}\ I\subset(0,1)\\ \iff
\\ &\ \tag{\Radioactivity}\sum_{k=1}^n\mu([\v_k=n+I])\sim u(n)|I|.
\end{align*}\endproclaim
These follow from the  ''remarks about mixing`` on pages \pageref{HKmixing} and \pageref{sfHKmixing}.
\

As shown in \cite{G-L}\ $\&$\ \cite{C-D}, (\faBolt)/(\Radioactivity) always hold when
$a(n)\gg\sqrt n$. Examples where they fail exist whenever $a(n)\ngg\sqrt n$.

In the case where $(\Om,\mu,\tau,\v)$ is the excursion time process from $0$ generated by an aperiodic, 
$\Bbb Z$-valued random walk in the domain of attraction of a symmetric, $p$-stable law ($1<p\le 2$),
the convergence version of (\ref{faLightbulbO}) follows from results in \cite{Wendel}, \cite{Liggett1970}, \cite{Vervaat}.
\subsection*{Examples:\ Intermittent interval maps}\ \ Consider the following 
piecewise onto\ 
interval maps $T_\g:[0,1]\to [0,1]$, ($\g>0$)\ \ defined by
\begin{align*}
		T_\g x:=
	\begin{cases}
		&x(1+(2x)^{1/\gamma}), \ \ \ 0\leq x<\frac{1}{2},
		\\ &
		2x-1,\ \ \ \ \frac{1}{2}\leq x\leq 1.
	\end{cases}
\end{align*}
Each $T_\g$ is conservative and exact (with respect to Lebesgue measure) and admits an invariant density $h_\g$ (unique up to constant multiplication) which is continuous on $(0,1)$ and satisfies  $h_\g(x)\propto x^{-1/\gamma}$ as $x\to 0$ (see \cite{Tha83}).
\

For $\g>1,\ h_\g$ is integrable and the maps $T_\g:\ \g>1$ were studied in \cite{LSV} as ''probabilistic intermittency``.
\

For $0<\g\le 1$,  the invariant measure $d\mu_\g(x)=h_\g(x)dx$ is infinite  and
(see e.g.  \cite{RFEx},\cite[\S4]{IET}):
\sbul\  $([0,1],\mu_\g,T_\g)$ is pointwise dual ergodic with $a_n(T_\g)\propto n^\g$ when $0<\g<1$ and
$a_n(T_1)\propto\tfrac{n}{\log n}$;
\sbul \  $T_{[\tfrac12,1]}$ is a  piecewise onto {\tt AFU} map (as on p. \pageref{AFU})\  with the return time function constant of intervals of monotonicity.
\

By Theorem B, $[\tfrac12,1]$ is a {\tt GLL} set for $T_\g$ and \ref{dsmilitary} holds by Theorem A.
\

Now let $\kappa\ge 2$ and consider the maps $R_\g:\ 0<\g\le 1$ defined by
\begin{align*}
		R_\g x:=
	\begin{cases}
		&x(1+(\kappa x)^{1/\gamma})\ \mod 1, \ \ \ 0\leq x<\frac{1}{2},
		\\ &
		2x-1,\ \ \ \ \frac{1}{2}\leq x\leq 1.
	\end{cases}
\end{align*}
These are (possibly non-Markov) 
 {\tt AFN} maps as in \cite{RolandAFN} and, as shown there are 
 absolutely continuous, $\s$-finite. $R_\g$-invariant measures $\nu_\g$ on $[0,1]$ ($0<\g\le 1$). Moreover, each
\sbul\  $([0,1],\nu_\g,R_\g)$ is conservative,  exact and pointwise dual ergodic with $a_n(R_\g)\propto n^\g$ when $0<\g<1$ and
$a_n(R_1)\propto\tfrac{n}{\log n}$.
\

It is also shown in \cite{RolandAFN} that a conservative {\tt AFN} map induces an {\tt AFU} map on some interval,
which in turn induces an exact {\tt AFU} map on some sub-interval (see p. \pageref{AFU}). This yields  {\tt GLL} sets
for each $R_\g$ and    establishes \ref{dsmilitary} as above.

\subsection*{Examples:\ Geodesic flows}\ \ In \S6, we consider an integrated, tied down mixing property \ref{faTrain} (on p. \pageref{faTrain}) which is satisfied by the natural extensions of transformations satisfying \ref{dsmathematical} (on p.
\pageref{dsmilitary}). In particular, (on p. \pageref{cyclic}),  we show that the geodesic flow (on the unit tangent bundle) of a cyclic cover of a
compact hyperbolic surface satisfies \ref{faTrain}.

\

\section*{\S2 Proof of Theorem A}
\

\f{\bf Opening Remark}\label{WLOG}\ \   We first note that to establish the conditional, tied-down renewal mixing property  on p. \pageref{dsmilitary}, it suffices to consider  (\ref{dsmilitary}) for some countable sub-collection $\mathcal A$ of $C([0,\infty])$ with dense linear span.
\

To see this, fix $C\in\mathcal F_+\ \&\ f\in L^1_+$ and suppose that that (\ref{dsmilitary})   holds for fixed $\forall\ g\in\mathcal A$. 
For $x\in X$ define the measures $P_{n,x}$ on $\Bbb R_+$ by 
$$P_{n,x}(g)=\int_XgdP_{n,x}:=\tfrac1{m(C)u(n)}\T^n(1_Cg(\tfrac{S_n(f)}{a(n)}))(x).$$
\

By \cite[Prop.4.2]{RatWM}, for each $g\in\mathcal A,\ \exists\ X_g\in\B(X),\ m(X\setminus X_g)=0$ so that for  $x\in X_g,\ \exists\ K_{g,x}\subset\Bbb N$ so that $\sum_{k=1}^n1_{K_{g,x}}(k)u(k)=o(a(n))$ and so that
$$P_{n,x}(g)\xrightarrow[n\to\infty,\ n\notin K_{g,x}]{}\Bbb \Bbb E(g(m(f)W_\g)).$$
    Set $X_\mathcal A:=\bigcap_{g\in\mathcal A}X_g$. 
    \
    
    Since $\mathcal A$ is countable $m(X\setminus X_\mathcal A)=0$. 
    \
    
    Next, for fixed $x\in X_\A$, we ''uniformize`` the sets $\{K_{g,x}:\ g\in\A\}$.
    \
    
    By \cite[Rem.4.1(iii)]{RatWM}, $\exists\ K_{x}\subset\Bbb N$ so that $\sum_{k=1}^n1_{K_{x}}(k)u(k)=o(a(n))$ and so that
    $\forall\ g\in\A,\ K_{g,x}\cap [n,\infty)\subset K_x$ for large $n\in\Bbb N$. 
    \
    
    It follows that
    $$P_{n,x}(g)\xrightarrow[n\to\infty,\ n\notin K_{x}]{}\Bbb \Bbb E(g(m(f)W_\g))\ \forall\ g\in\A,$$
    whence by uniform approximation $\forall\ g\in C([0,\infty])$, and by monotone approximation in $L^1(\text{\tt dist.}\,W_\g),\ \forall\ g\in C_B(\Bbb R_+)$.
    \
    
    Thus, again by \cite[Prop. 4.2]{RatWM},   (\ref{dsmilitary}) holds $\forall\ g\in C_B(\Bbb R_+)$.

 \

Let $\Om\in \mathcal F$ be a {\tt GLL} set with accompanying $T_\Om$-generating partition $\a$.  
Up to isomorphism,
$\Om=\a^\mathbb N$,  $T_\Om:\Om\to \Om$
is the shift and the collection $\mathcal C_\a(T_\Om)$ of $(\a,T_\Om)$-cylinder sets forms a base of clopen sets for the Polish topology on $\Om$.
\

\proclaim{Lemma 2.1}\label{lemma2.1}
\

Let the {\tt CEMPT}  $(X,m,T)$ be pointwise dual ergodic  with $\g$-regularly varying return sequence   
$a(n)=a_n(T)$   {\rm($\g\in (0,1)$)}. 
\

Suppose that $(X,m,T)$  has a  {\tt GLL} set 
 $\Om\in \B(X),\ 0<m(\Om)<\infty$, then for $A\in\mathcal C_\a(T_\Om),\ g\in C_B(\Bbb R),\ g\ge 0$,
\begin{align*}&\tag{GL}\varliminf_{n\to\infty}\tfrac1{u(n)}\T^n(1_Ag(\tfrac{S_n(1_\Om)}{a(n)}))\ge m(A) \Bbb E(g(W_\g))\ \text{a.e. on $\Om$,}
\\ &\tag{\Wheelchair}\tfrac1{a(N)}\sum_{n=1}^N\T^n(1_Ag(\tfrac{S_n(1_\Om)}{a(n)}))
\xrightarrow[N\to\infty]{}\  m(A) \Bbb E(g(W_\g))\ \text{a.e. on $\Om$,}
\end{align*} where  $\a$ is the accompanying $T_\Om$-generating partition and $u(n)\sim\tfrac{\g a(n)}n$.\endproclaim
 \demo{Proof}
 \
 
 We begin with  {\rm (GL)} and consider only the periodic case, the aperiodic case being similar and easier 
 (c.f. Lemma 2.2.1 in \cite{G-L} and Lemma 9.2 in \cite{RatWM}).
 \

\

   Fix $A\in\mathcal C_\a(T_\Om)$ and $0<c<d<\infty$ and  $g\in C_B(\Bbb R_+)_+$ $[c,d]$-{\it nice}\label{cdnice} in the sense that  $\log G:[c,d]\to\Bbb R$ is smooth where
$G(x):=g(\tfrac1{x^\g})f_{Z_\g}(x)$ with $f_{Z_\g}$ the probability density of $Z_\g$.
\

Writing $x_{k,n}:=\tfrac{n}{a^{-1}(k)}$ for $1\le k\le n$, we have by 
 regular variation that
 $$k\sim \frac{a(n)}{x_{k,n}^\g}\ \text{as}\ \ k,\ n\to\infty,\ x_{k,n}\in [c,d].$$
\

Using this and the {\tt GLL} property of $\Om$, we have, as $n\to\infty$,
\begin{align*}\tag{{\scriptsize\faTaxi}}\label{faTaxi}\T^n(1_Ag(\tfrac{S_n(1_\Om)}{a(n)}))&=\sum_{k=1}^n\T_\Om^{k}(1_{A\cap[\v_k=n]}g(\tfrac{k}{a(n)}))\\ &\ge
\sum_{1\le k\le n,\ x_{k,n}\in [c,d]}\T_\Om^{k}(1_{A\cap[\v_k=n]}g(\tfrac{k}{a(n)}))\\ &\sim
\sum_{1\le k\le n,\ x_{k,n}\in [c,d]}g(\tfrac1{x_{k,n}^\g})\T_\Om^{k}(1_{A\cap[\v_k=n]})\\ &\sim
m(A)\sum_{1\le k\le n,\ x_{k,n}\in [c,d]}g(\tfrac1{x_{k,n}^\g})\tfrac{pf_{Z_\g}(x_{k,n})}{a^{-1}(k)}1_{p\Bbb Z}(n-k\xi).
\end{align*}
To continue, we'll need
\proclaim{Lemma 2.2}\ 
\

Suppose that $a(n)$ is $\g$-regularly varying with $\g\in (0,1)$ and satisfies 
\begin{align*}
 \tag{{\scriptsize\faPlug}}\label{faPlug} a^{-1}(n+1)-a^{-1}(n)\sim \tfrac{a^{-1}(n)}{\g n}.
\end{align*}
Let $\Bbb G=\Bbb Z$ or $\Bbb R$, and let $\xi,\ p\in\Bbb G.\ 0<\xi<p$ satisfy \f$\overline{\<\{n\xi:\ n\ge 1\}+p\Bbb Z\>}=\Bbb G$, then,
for   $0<c<d<\infty$,
 $g\in C_B(\Bbb R_+)$ $[c,d]$-nice  and 
 $I\subset (0,p)\cap\Bbb G$ an interval:
\begin{align*}\tag{\dsarchitectural}\label{dsarchitectural}\frac1{u(n)}\sum_{1\le k\le n,\ x_{k,n}\in [c,d]}&g(\tfrac1{x_{k,n}^\g})\tfrac{pf_{Z_\g}(x_{k,n})}{a^{-1}(k)}1_{I+p\Bbb Z}(n-k\xi)
 \\ &\xrightarrow[n\to\infty]{}\ \mathbb m_{\Bbb G}(I)\Bbb E(1_{[c,d]}(Z_\g)g(Z_\g^{-\g})Z_\g^{-\g}).
\end{align*} where $u(n):=\tfrac{\g a(n)}n$ and $x_{k,n}:=\tfrac{n}{a^{-1}(k)}$ for $1\le k\le n$.\endproclaim
Note that \ref{faPlug} can be ensured   by possibly passing to an asymptotically equivalent function.
\

We'll prove Lemma 2.2 after the proof of Lemma 2.1.
\

\demo{Proof of {\rm (GL)} given Lemma 2.2}
By Lemma 2.2, for $g\in C_B(\Bbb R_+)$ $[c,d]$-nice (as on p. \pageref{cdnice}):
\begin{align*}\frac1{u(n)}\sum_{1\le k\le n,\ x_{k,n}\in [c,d]}&g(\tfrac1{x_{k,n}^\g})\tfrac{pf_{Z_\g}(x_{k,n})}{a^{-1}(k)}1_{p\Bbb Z}(n-k\xi)
 \\ &\xrightarrow[n\to\infty]{}\  \Bbb E(1_{[c,d]}(Z_\g)g(Z_\g^{-\g})Z_\g^{-\g}).
\end{align*}

Now $$ \Bbb E(1_{[c,d]}(Z_\g)g(Z_\g^{-\g})Z_\g^{-\g})
\xrightarrow[c\to 0+,\ d\to\infty]{}\   \Bbb E(g(Y_\g )Y_\g)= \Bbb E(g(W_\g)),$$
$$\therefore\ \T^{n}(1_Ag(\tfrac{S_n(1_\Om)}{a(n)}))\gtrsim\ \ u(n)m(A) \Bbb E(g(W_\g)).\ \ \CheckedBox\text{\rm (GL)}$$
\demo{Proof of (\Wheelchair) given {\rm(GL)}}\ \ It suffices to show that a.e.,
\begin{align*}\varlimsup_{N\to\infty}\tfrac1{a(N)}\sum_{n=1}^N\sum_{1\le k\le n,\ x_{k,n}\notin [c,d]}\T_\Om^{k}(1_{A\cap[\v_k=n]}g(\tfrac{k}{a(n)}))\xrightarrow[c\to 0,\ d\to\infty]{}\ 0. 
\end{align*}
To this end, note that
$$\sum_{1\le k\le n,\ x_{k,n}\notin [c,d]}\T_\Om^{k}(1_{A\cap[\v_k=n]}g(\tfrac{k}{a(n)}))\le
\|g\|_\infty\sum_{1\le k\le n,\ x_{k,n}\notin [c,d]}\T_\Om^{k}(1_{A\cap[\v_k=n]})$$
and
\begin{align*}\sum_{n=1}^N&\sum_{1\le k\le n,\ x_{k,n}\notin [c,d]}\T_\Om^{k}(1_{A\cap[\v_k=n]})\\ &=
 \sum_{n=1}^N\sum_{1\le k\le n}\T_\Om^{k}(1_{A\cap[\v_k=n]})-\sum_{n=1}^N\sum_{1\le k\le n,\ x_{k,n}\in [c,d]}\T_\Om^{k}(1_{A\cap[\v_k=n]})\\ &=\sum_{n=1}^N\T^n1_A-\sum_{n=1}^N\sum_{1\le k\le n,\ x_{k,n}\in [c,d]}\T_\Om^{k}(1_{A\cap[\v_k=n]})\\ &=\D_N(c,d).
\end{align*}
Now  by (GL) with $g\equiv 1$,
$$\frac{\D_N(c,d)}{a(N)}\lesssim\ m(A)(1- \Bbb E(1_{[c,d]}(Z_\g)Z_\g^{-\g}))
\xrightarrow[c\to 0,\ d\to\infty]{}\ 0.\ \CheckedBox\ \text{(\Wheelchair)}$$
\

The proof of Lemma 2.2  makes use of  the following standard
\proclaim{Equidistribution Lemma}
\

Suppose that $K$ is a compact, Abelian group and that \ $\xi\in K,\ \overline{\{n\xi:\ n\ge 1\}}=K$.
\

If $u_n^{(\nu)}\ge 0,\ \ (n,\ \nu\ge 1)$ satisfies
\sms {\rm (i)} $\sum_{n\ge 1}u_n^{(\nu)}\xrightarrow[\nu\to\infty]{}\ C\in \Bbb R_+\ \&$;  {\rm (ii)}  $\sum_{n\ge 1}|u_n^{(\nu)}-u_{n+1}^{(\nu)}|\xrightarrow[\nu\to\infty]{}0$,
\

then 
\begin{align*}\tag{\Info}\label{Info}\sum_{n\ge 1}u_n^{(\nu)}1_U(x+n\xi)&\xrightarrow[\nu\to\infty]{}Cm_K(U)
\ \ \forall\  U\in\B(K),\ m_K(\bdy U)=0 \\ & \text{uniformly in}\ \ x\in K
\end{align*}where $m_K$ denotes normalized Haar measure on $K$.\endproclaim\demo{Proof sketch}\ \ 
\

Define $r:K\to K$ by $r(x)=x+\xi$. Since $\overline{\{n\xi:\ n\ge 1\}}=K$, the only $r$-invariant probability on $K$ is the normalized Haar measure $m_K$.
\

Define probabilities $p_{\nu,x}\in\mathcal P(K)$ ($x\in K,\ \nu\ge 1$) by
$$p_{\nu,x}(U):=\ 
\tfrac1{C_\nu}\sum_{n\ge 1}u_n^{(\nu)}1_U(x+n\xi)$$ where $C_\nu:=\sum_{n\ge 1}u_n^{(\nu)}$.
\

If\  (\Info)\  fails, $\exists\ \e>0,\ x_k\in X,\ n_k\to\infty\ \&\ f\in C(K)$ so that 
$$|p_{n_k,x_k}(f)-m_K(f)|\ge\e\ \forall\ k\ge 1.$$
By compactness, by passing to a subsequence, it can be ensured that 
$$x_k\xrightarrow[k\to\infty]{}x\in K\ \&\ p_{n_k,x_k}\xrightarrow[k\to\infty]{}P\in\mathcal{P}(K)$$ weakly, whence
$|P(f)-m_K(f)|\ge\e$.
\

On the other hand, it follows from  (ii) that 
$$p_{n_k,x_k}(g)-p_{n_k,x_k}(g\circ r)\xrightarrow[k\to\infty]{}0\ \forall\ g\in C(K)$$
whence $P=m_K$ contradicting failure of (\Info).\ \ \Checkedbox\ 
\

\subsection*{Proof of Lemma 2.2}
\

We  establish \ref{dsarchitectural} (p. \pageref{dsarchitectural}).
\

By \ref{faPlug}, 
$$x_{k,n}-x_{k+1,n}=\frac{n}{a^{-1}(k)}-\frac{n}{a^{-1}(k+1)}\sim \frac{n}{\g ka^{-1}(k)}$$
as $k,\ n\to\infty$, $x_{k,n}\in [c,d]$.
\

Also
$$a(n)=a(x_{k,n}a^{-1}(k))\sim x_{k,n}^\g a(a^{-1}(k))\sim x_{k,n}^\g k$$
as $k,\ n\to\infty,\ x_{k,n}\in [c,d]$ by the uniform convergence theorem for regularly varying functions.
Thus:
\begin{align*}\tag{{\scriptsize\faLinux}}\label{faLinux}
 \tfrac1{a^{-1}(k)}\sim\tfrac{\g k}n\cdot (x_{k,n}-x_{k+1,n})\sim \tfrac{\g a(n)}n\cdot\tfrac{x_{k,n}-x_{k+1,n}}{x_{k,n}^\g}=u(n)\cdot\tfrac{x_{k,n}-x_{k+1,n}}{x_{k,n}^\g}
\end{align*}

whence, as $n\to\infty$,
                        
\begin{align*}\sum_{1\le k\le n,\ x_{k,n}\in [c,d]}&g(\tfrac1{x_{k,n}^\g})\frac{pf_{Z_\g}(x_{k,n})}{a^{-1}(k)}1_{I+p\Bbb Z}(n-k\xi)\\ &\sim
pu(n)
\sum_{1\le k\le n,\ x_{k,n}\in [c,d]}
g(\tfrac1{x_{k,n}^\g})\tfrac{(x_{k,n}-x_{k+1,n})}{x_{k,n}^\g}f_{Z_\g}(x_{k,n})1_{I+p\Bbb Z}(n-k\xi).
\end{align*}
In order to use the Equidistribution Lemma(p. \pageref{Info}), define
\begin{align*}
 v^{(n)}_k&:=g(\tfrac1{x_{k,n}^\g})\tfrac{pf_{Z_\g}(x_{k,n})}{a^{-1}(k)}1_{[x_{k,n}\in [c,d]]}=
 \tfrac{pH(x_{k,n})}n1_{[x_{k,n}\in [c,d]]}\\ &\text{with}\ \ H(x):=xg(\tfrac1{x^\g})f_{Z_\g}(x).
\end{align*}

By (\ref{faLinux}), 
as $k,\ n\to\infty,\ x_{k,n}\in [c,d]$,
\begin{align*}
 v^{(n)}_k\sim u(n)g(\tfrac1{x_{k,n}^\g})\tfrac{pf_{Z_\g}(x_{k,n})(x_{k,n}-x_{k+1,n})}{x_{k,n}^\g}1_{[x_{k,n}\in [c,d]]}.
\end{align*}
Thus, using the convergence of Riemann sums to the  Riemann integral, 
\begin{align*}\tag*{{\scriptsize\faPaw}}\label{faPaw}
 \sum_{k\ge 1} v^{(n)}_k\sim pu(n) \Bbb E(1_{[c,d]}(Z_\g)g(Z_\g^{-\g})Z_\g^{-\g}).
\end{align*}

Next suppose that $x_{k_0-1,n}<c\le x_{k_0,n}<x_{k_1,n}\le d<x_{k_1+1,n}$.
For $k\in [k_0,k_1)$, 
\begin{align*}|v_k^{(n)}-v_{k+1}^{(n)}|&=\tfrac{p}n |H(x_{k,n})-H(x_{k+1,n})|\\ &\sim
\tfrac{p}n |H'(x_{k,n})|(x_{k,n}-x_{k+1,n})\\ &\le M\tfrac{p}n (x_{k,n}-x_{k+1,n})
\end{align*}
where $M:=\sup_{x\in [c,d]}|H'(x)|$
and
\begin{align*}
 \sum_{k\ge 1}|v_k^{(n)}-v_{k+1}^{(n)}|&\le \tfrac{pM}n\sum_{k=k_0}^{k_1-1}(x_{k,n}-x_{k+1,n})+|v_{k_0}^{(n)}|+|v_{k_1}^{(n)}|\\ &=
 \tfrac{pM}n(x_{k_0,n}-x_{k_1,n})+\tfrac{2pK}n\ \text{with} \ K:=\sup_{x\in [c,d]}|H(x)|\\ &\le
 \frac{R}n\ \text{with}\ R:=p(M(d-c)+2K)\\ &=o(u(n))\ \text{as}\ n\to\infty.
\end{align*}

Thus, by  (\Info) with $u_k^{(n)}:=\frac{v_k^{(n)}}{u(n)}$,
                                    
\begin{align*}
\sum_{k\ge 1,\ x_{k,n}\in [c,d]}&g(\tfrac1{x_{k,n}^\g})\tfrac{pf_{Z_\g}(x_{k,n})}{a^{-1}(n)}
1_{I+p\Bbb Z}(n-k\xi)=\sum_{k\ge 1}v_k^{(n)}1_{I+p\Bbb Z}(n-k\xi)\\ &
\sim pu(n) \Bbb E(1_{[c,d]}(Z_\g)g(Z_\g^{-\g})Z_\g^{-\g})\ \cdot\ m_{\Bbb G/p\Bbb Z}(I+p\Bbb Z)\\ &=
u(n)m_\Bbb G(I) \Bbb E(1_{[c,d]}(Z_\g)g(Z_\g^{-\g})Z_\g^{-\g}).\ \ \CheckedBox\ \text{(\dsarchitectural)}\end{align*}
\f{\bf Remark}\ The  aperiodic form of Lemma 2.2 is \ref{faPaw}\ with $p=1$.

\

\subsection*{ Remarks about mixing}\label{HKmixing} \ \ \
\

It follows from (GL) $\&$ (\ref{faTaxi}) on p. \pageref{faTaxi} that for $(X,m,T)$  pointwise, dual ergodic {\tt MPT} with $a(n)=a_n(T)$ $\g$-regularly varying ($\g\in (0,1)$), and  a 
{\tt GLLT} set $\Om\in\B(X),\ 0<m(\Om)<\infty$,  TFAE:
\begin{align*}&\tag{i}\tfrac1{u(n)}\T^n(1_Ag(\tfrac{S_n(1_\Om)}{a(n)}))\xrightarrow[n\to\infty]{} m(A) \Bbb E(g(W_\g))\ \ \text{ a.e. on $\Om$}
\\ & \ \ \ \ \ \ \ \ \ \ \ \ \ \ \ \ \ \ \forall\ A\in\mathcal C_\a(T_\Om),\ g\in C_B(\Bbb R),\ g\ge 0;\\ &\tag{ii}\varlimsup_{n\to\infty}\tfrac1{u(n)}\sum_{1\le k\le n,\ x_{k,n}\notin [c,d]}\T_\Om^{k}(1_{A\cap[\v_k=n]})\xrightarrow[c\to 0+,\ d\to\infty]{}\ 0\ \text{a.s. on}\ \Om;\\ &\tag{iii} \T^n1_\Om\sim u(n)\ \ \ \text{a.s. on}\ \Om.
\end{align*}
   
By \cite[theorem 2.1]{Melter}, (iii) holds (whence also (i))  if   $\g\in (\tfrac12,1)$, $T_\Om$ is mixing and $\exists$ a partition $\a\subset\B(\Om)$ so that
 $(\Om,m_\Om,T_\Om,\a)$ is  either a Gibbs-Markov or an AFU map with $\varphi:\Om\to\Bbb N$ non-arithmetic and $\a$-measurable.

\demo{Proof of Theorem A:  finish}\label{A-finish} 
\

It suffices to prove (\ref{dsmilitary}) for  $g>0$ smooth, with
 $x\mapsto\ \log g(e^x)$  uniformly continuous on $[-\infty,\infty]$ (see the opening remark on  p. \pageref{WLOG}).
 \
 
 \Par\  We first fix  $f=1_\Om$ and $g>0$ as above; and show that (\ref{dsmilitary}) (on p. \pageref{dsmilitary}) holds for  $C\in\mathcal F_+$.
 \
 
 To this end, let
 $$\frak R=\frak R(g):=\{A\in\mathcal F_+:\ \text{(\ref{dsmilitary}) holds for}\ 1_\Om,\ g\ \&\ A\}.$$
 We'll show that  $\frak R=\mathcal F_+$, and first show that $\frak R$ contains a dense, hereditary ring. 
 \
 
By Lemma 2.1,  we see using  of \cite[Prop. 4.2]{RatWM} as in the proof of \cite[Prop. 3.1]{nice},  that
$\mathcal B(\Om)\subset\frak R$. 
\

Next, we claim that  $T^{-1}\frak R\subset\frak R$. 
\

To see this, let $C\in\frak R$, then
\begin{align*}
\T^n(1_{T^{-1}C}g(\tfrac{S_n(1_\Om)}{a(n)}))&=\T^n(1_C\circ Tg(\tfrac{1_\Om+S_{n-1}(1_\Om)\circ T}{a(n)}))\\ &\sim \T^n(1_C\circ Tg(\tfrac{S_{n-1}(1_\Om)\circ T}{a(n-1)}))
\ \text{by u.c. of } \ x\mapsto\ \log g(e^x)\\ &=\T^{n-1}(1_Cg(\tfrac{S_{n-1}(1_\Om)}{a(n-1)}))
\end{align*}
and $T^{-1}C\in\frak R$.
\

Thus $\frak R$ contains the dense, hereditary ring 
$$\mathcal R:=\{A\in\mathcal{F}_+:\ \exists\ n\ge 0,\ A\subset \bigcup_{k=0}^nT^{-k}\Om\}.$$
Now fix  $C\in\mathcal F_+$ and let  $D\in\mathcal R$ so that
$D\subset C$.
\

It follows that
\begin{align*}|\T^n(1_C&g(\tfrac{S_n(1_\Om)}{a(n)}))-u(n)m(C)\Bbb E(g(W_\g))|\\ &\le
\T^n(1_{C\setminus D}g(\tfrac{S_n(1_\Om)}{a(n)}))+|\T^n(1_Dg(\tfrac{S_n(1_\Om)}{a(n)}))-u(n)m(D)\Bbb E(g(W_\g))|\\ &
\ \ \ \ \ \ \ \ \ \ \ +u(n)m(C\setminus D)\Bbb E(g(W_\g)) \\ &\le
|\T^n(1_Dg(\tfrac{S_n(1_\Om)}{a(n)}))-u(n)m(D)\Bbb E(g(W_\g))|\\ &
\ \ \ \ \ \ \ \ \ \ \ +\|g\|_{C_B}(\T^n1_{C\setminus D}+u(n)m(C\setminus D)),
\end{align*}
whence using $D\in\mathcal R\subseteq\frak R$,
\begin{align*}\varlimsup_{N\to\infty}\tfrac1{a(N)}&\sum_{n=1}^N|\T^n(1_Cg(\tfrac{S_n(1_\Om)}{a(n)})-u(n)m(C)\Bbb E(g(W_\g))|\\ &
\le \|g\|_{C_B}\varlimsup_{N\to\infty}\tfrac1{a(N)}\sum_{n=1}^N(\T^n1_{C\setminus D}+u(n)m(C\setminus D)\\ &=
2\|g\|_{C_B}m(C\setminus D)\ \ \text{a.e.  by {\tt PDE} of $T$}\\ &\xrightarrow[D\uparrow C,\ D\in\mathcal R]{} 0
\end{align*}
and $\frak R=\mathcal F_+$.\ \Checkedbox\P
 
 \
 
To complete the proof of (\dsmilitary), fix $F\in L^1_+$ and $g>0$ smooth, with
 $x\mapsto\ \log g(e^x)$  uniformly continuous on $[-\infty,\infty]$. 
 \
 
 By the ratio ergodic theorem
$$\frac{S_n(F)}{S_n(1_\Om)}\xrightarrow[n\to\infty]{}\ m(F) \ \text{a.e.}$$
Suppose the convergence is uniform on $C\in\mathcal F_+$, then $\exists\ \e_N\downarrow 0$ so that
$$g(\tfrac{S_n(F)}{a(n)})=(1\pm\e_n)g(\tfrac{m(F)S_n(1_\Om)}{a(n)})\ \text{on}\ C\ \forall\ n\ge 1.$$
Thus
\begin{align*}|\T^n(1_Cg(&\tfrac{S_n(F)}{a(n)}))-m(C) \Bbb E(g(m(F)W_\g))u(n)|\\ & \le 
|\T^n(1_{C}(g(\tfrac{S_n(F)}{a(n)})-g(\tfrac{m(F)S_n(1_\Om)}{a(n)})))|\\ &\ \ \ \
+|\T^n(1_Cg(\tfrac{m(F)S_n(1_\Om)}{a(n)}))-m(C) \Bbb E(g(m(F)W_\g))u(n)|\\ &\le
\e_n\|g\|_\infty\T^n1_C+|\T^n(1_Cg(\tfrac{m(F)S_n(1_\Om)}{a(n)}))-m(C) \Bbb E(g(m(F)W_\g))u(n)|\end{align*}
and using \P,
\begin{align*}\tfrac1{a(N)}\sum_{n=1}^N&|\T^n(1_Cg(\tfrac{S_n(F)}{a(n)}))-m(C) \Bbb E(g(m(F)W_\g))u(n)|\\ &\le
\tfrac1{a(N)}\sum_{n=1}^N(\e_n\|g\|_\infty\T^n1_C+|\T^n(1_Cg(\tfrac{m(F)S_n(1_\Om)}{a(n)}))-m(C) \Bbb E(g(m(F)W_\g))u(n)|)\\ &
\xrightarrow[N\to\infty]{}\ 0\ \text{a.e.}.
\end{align*}
Now fix $D\in\mathcal F_+$. By Egorov's theorem $\exists\ C_\nu\in\mathcal F_+,\ C_\nu\uparrow D\ \mod\ m$ so that
$$\frac{S_n(F)}{S_n(1_\Om)}\xrightarrow[n\to\infty]{}\ m(F) \ \text{uniformly on each}\ C_\nu.$$
Thus,  a.e., $\forall\ \nu\ge 1$,  as $N\to\infty$,
\begin{align*}\tfrac1{a(N)}\sum_{n=1}^N&|\T^n(1_Dg(\tfrac{S_n(F)}{a(n)}))-m(D) \Bbb E(g(m(F)W_\g))u(n)|\\ &\lesssim
\tfrac{\|g\|_\infty}{a(N)}\sum_{n=1}^N(\T^n1_{D\setminus C_\nu}+m(D\setminus C_\nu)u(n))\\ &
\xrightarrow[N\to\infty]{}\ 2\|g\|_\infty m(D\setminus C_\nu)\\ & \xrightarrow[\nu\to\infty]{}\ 0.\ \text{\Checkedbox (\dsmilitary)}
\end{align*}

\section*{\S3 Periodic {\tt LLT}s}
In this section, we give conditions for an asymptotically $\g$-stable, positive, ergodic, stationary  process $(\Om,\mu,\tau,\phi)$
to satisfy a periodic {\tt LLT} as defined on p. \pageref{LLT}. 

\subsection*{Fibered systems}\label{fibered}
\

We assume first that there is a  countable, generating partition $\alpha\subset\B(\Om)$ so that
$\tau:a\to\tau a$ is  invertible and nonsingular for $a\in\alpha$. Such  $(\Om,\mu,\tau,\a)$ is called a {\it fibered system}.
\

Up to isomorphism,  a fibered system  $(\Om,\mu,\tau,\a)$ has the form:
\sms $\Om\subset S^\Bbb N$ where is a $S$ is a 
finite or countable 
{\it state space} and $\Om$ is a {\it subshift} (i.e.  shift invariant and closed with respect to the 
 polish product topology on  $S^\Bbb N$) is a  closed,  subset of  subshift 
 and
$$\a=\{[s]=\{(x_1,x_2,\dots)\in X:\ x_1=s\}:\ s\in S\}.$$

Thus the fibered system $(\Om,\mu,\tau,\a)$ can be considered as a continuous map of a polish space.

If $(\Om,\mu,\tau,\a)$ is a fibered system, then for $n\ge 1$, so is 
\f $(\Om,\mu,\tau^n,\alpha_n)$ where $\alpha_n:=\bigvee_{k=0}^{n-1}\tau^{-k}\alpha$.
\subsection*{Inverse branches and the transfer operator}

For $n\ge 1$, there are $\mu$-nonsingular inverse branches of $\tau^n$
denoted $v_a:\tau^na\to a\ (a\in\a_n)$ with Radon Nikodym derivatives
$$v_a':={d\mu\circ v_a\over d\mu}:\tau^na\to\Bbb R_+.$$
The transfer operator is given by
$$\ttau(f):=\sum_{a\in\a}1_{\tau a}v_a'f\circ v_a.$$

\subsection*{Doeblin-Fortet operators}\label{quasicompact}\footnote{see \cite{AD} and references therein, also \cite{Norman,HH}}
\

Let $\mathcal L\subset L^1(\mu)$ be a Banach space  so that 
$(L^1(\mu),\mathcal{L})$ is an 
{\it adapted pair} in the sense that
$\mathcal L\subset {L^1(\mu)},\ \|\cdot\|_{{L^1(\mu)}}\le\|\cdot\|_{\mathcal L},\ (\o {\mathcal L})_{{L^1(\mu)}}={L^1(\mu)}$,
and $\mathcal L$-closed,  $\mathcal L$-bounded sets are ${L^1(\mu)}$-compact.
\

 We say that an operator $P\in\hom(\mathcal L,\mathcal L)\cap\hom( L^1(\mu),L^1(\mu)) $ is {\it Doeblin-Fortet} (D-F) on $(L^1(\mu),\mathcal{L})$  if
 \begin{align*}\tag*{{\tt DF}(i)}& \|P^nf\|_1\le H\|f\|_1\ \forall n\in\Bbb N,\ f\in L^1(\mu)\ \&\\ &\tag*{{\tt DF}(ii)}
\|Pf\|_{\mathcal L}\le \th \|f\|_{\mathcal L}+R\|f\|_1\ \forall\ f\in \mathcal L.  
 \end{align*}

where $R,H\in\Bbb R_+$ and $\th\in (0,1)$.
\

It shown in \cite{IT-M} that a D-F operator $P\in\hom(\mathcal L,\mathcal L)$  has spectral radius $\rho(P)\le 1$ and that,
if $\rho(P)=1$, then $P$ is {\it quasicompact} in the sense that
$\exists$ $A=A(P)\in\hom(\mathcal L,\mathcal L)$ of form
$$A=\sum_{k=1}^N\l_kE_k$$ with $N\ge 1,\ E_1,\dots,E_N\in\hom(\mathcal L,\mathcal L)$ finite dimensional  projections,
$\l_1,\dots,\l_N\in S^1:=\{z\in\Bbb C:\ |z|=1\}$  
so that $\rho(P-A)<1$.
\

In particular for $\rho(P-A)<\th<1,\ \exists\ M>0$ so that 
$$\|P^nf-A^nf\|_\mathcal{L}\le M\th^n\|f\|_\mathcal{L}\ \ \forall\ n\ge 1,\ f\in \mathcal{L}.$$
Moreover, the spectral radius of a D-F operator $P\in\hom(\mathcal L,\mathcal L)$ satisfies $\rho(P)\le 1$ with equality iff $\exists\ \l\in\Bbb S^1\ \&\ g\in L^1(\mu)$ so that
$Pg=\l g$, in which case $|g|$ is constant, $g\in\mathcal{L}$ and
$\|g\|_\mathcal L\le \tfrac{RH}{1-\th}\|g\|_1$ with $\th,R,H$ as in
{\tt DF}(ii).

\

For $P\in\hom(\mathcal L,\mathcal L)$ D-F with $\rho(P)=1$ we'll write  $\dim\,P:=\dim\,A(P)\mathcal{L}$.
\

If $\dim\,P=1$, then $A(P)=\l(P)N(P)$ where $\l(P)\in\Bbb C,\ |\l(P)|\le 1$ and $N(P)\in\hom(\mathcal L,\mathcal L)$ is a projection onto a  one dimensional subspace.
Let
$$\text{\tt DF}(1):=\{P\in\hom(\mathcal L,\mathcal L):\ {\tt D-F},\ \rho(P)=1\ \&\ \dim\,P=1\}.$$

The mixing of $\tau$ ensures that its transfer operator $\ttau$ is D-F, then $\ttau\in\text{\tt DF}(1)$ with $\l=1\ \&\ N(\ttau)f=\int_\Om fd\mu$.
We'll be interested in fibered systems for which the transfer operator is {\tt DF} on some  $(L^1(\mu),\mathcal{L})$.

\subsection*{Doeblin-Fortet fibered systems}\label{DFfibered}
\

We'll call the probability preserving, mixing   fibered system $(\Om,\mu,\tau,\a)$ a {\it Doeblin-Fortet {\rm ({\tt D-F})} fibered system} if there is 
a Banach space $L_\tau$ so that $(L^1(\mu),L_\tau)$  is an adapted pair, $1_A\in L_\tau\ \forall\ A\in\a$ and
the transfer operator $\ttau$ is {\tt DF} on $(L^1(\mu),L_\tau)$.
\subsection*{Example 1:\ {\tt AFU} maps}\label{AFU}
\

As in \cite{RolandNonLin}, an {\it {\tt AFU} map} is an  interval map   $(\Om,\mu,\tau,\a)$ where  for each $a\in\alpha$, $\tau|_a$ is (the restriction of) a  $C^2$ diffeomorphism 
$\tau:\overline{a}\to \overline{\tau a}$
 satisfying in addition:
                
\begin{align*}&\tag{A}\sup_X\tfrac{|\tau''|}{(\tau')^2}<\infty,\\ &
 \tag{F}\tau\a:=\{\tau A:A\in\alpha\}\ \text{\rm is finite},\\ &\tag{U}\inf_{x\in a\in\alpha}|\tau'(x)|>1.
\end{align*}
If  $(\Om,\mu,\tau,\a)$ is an {\tt AFU} map, then so is $(\Om,\mu,\tau^N,\a_N)$ for $N\ge 2$ with
$\a_N:=\bigvee_{k=0}^{N-1}\tau^{-k}\a$.
\

It follows from  \cite{Rychlik, RolandNonLin}  that:
\

\sbul a totally ergodic {\tt AFU} map $(\Om,\mu,\tau,\a)$  is exact, and a {\tt D-F} fibered system  with $L_\tau =\text{\tt BV}$, the space of functions on the interval $\Om$ with bounded variation equipped with the norm  $\|f\|_{L_\tau}:=\|f\|_1+\bigvee_\Om f$;

\sbul any {\tt AFU} map induces an exact {\tt AFU} map on some interval.

\subsection*{Example 2:\ Gibbs-Markov maps}\label{GibbsMarkov}
\


  The fibered system $(\Om,\mu,\tau,\a)$ is a {\it Markov map} if
  $\tau a\in\s(\alpha)\ \mod \mu\ \ \forall\ a\in\a$, and 
a {\it Gibbs-Markov {\rm ({\tt G-M})} map} if, in addition, $\inf_{a\in\alpha}\mu(\tau a)>0$ 

 and, for some $\th\in (0,1)$
 \begin{align*}\sup_{n\ge 1,\ a\in\alpha_0^{n-1},\ x,y\in \tau^na}
 \frac1{\th^{t(x,y)}}\cdot \left|\frac{v_a'(x)}{v_a'(y)}-1\right|<\infty
\end{align*}
with $t(x,y)=\min\{n\ge 1:\ \a_n(x)\ne \a_n(y)\}\le\infty$.
\

Note that a Markov {\tt AFU} map is Gibbs-Markov.
\

As shown in \cite[\S1]{AD}, a {\tt G-M} map $(\Om,\mu,\tau,\a)$  is a {\tt D-F} fibered system  with $L_\tau$, the space of $(\a,\th)$-{\it H\"older}\label{Holder} functions on $\Om$; that is
$$\{f:\Om\to\Bbb R:\ D_{\Om,\th,\a}(f):=\sup_{x,y\in \Om}\tfrac{|f(x)-f(y)|}{\th^{t(x,y)}}<\infty\}$$ 
equipped with the norm $\|f\|_{L_\tau}:=\|f\|_1+D_{\Om,\th,\a}(f)$.

\subsection*{Doeblin-Fortet cocycles}\label{DFcocycle} \ 
\

Let $(\Om,\mu,\tau,\a)$ be a {\tt D-F}  fibered system.
\

Let $\Bbb G\le\Bbb R^d$ be a subgroup of full dimension. For $\phi:\Om\to\Bbb G$ measurable and $\chi\in
\widehat{\Bbb G}$, define $P_{\phi,\chi}\in\hom(L^1(\mu),L^1(\mu))$ by 
$P_{\phi,\chi}(f):=\ttau(\chi(\phi)f)$.\label{pert}
\

Call the measurable function (cocycle) $\phi=(\phi^{(1)},\dots,\phi^{(d)}):\Om\to\Bbb G$
\bul {\it Doeblin-Fortet } ({\tt D-F}) if it is {\it locally $L_\tau$} in the sense that  $\exists\ M>0$ so that $$\forall\ 1\le k\le d\ \&\ A\in\a,\ 1_A\phi^{(k)}\in L_\tau\  \&\ 
\|1_A\phi^{(k)}-\int_A\phi^{(k)}d\mu\|_{L_\tau}\le M;$$
\

and
רעות 
\bul {\it admissible} if
\sms (i) $P_{\phi,\chi}\in\text{\tt DF}(1)\ \forall\ \chi\in\widehat{\Bbb G}$ and  \f(ii) $\chi\mapsto P_{\phi,\chi}$ is continuous ($\widehat{\Bbb G}\to \hom(L_\tau,L_\tau)$).
\

It follows from \cite{AD} and \cite{ADSZ}  (respectively) that a {\tt D-F} cocycle over a {\tt G-M} --  /  {\tt AFU}--  map is admissible.
\

\subsection*{Discrete {\tt D-F} cocycles}
\

For $(\Om,\mu,\tau,\a)$ a {\tt G-M} map, a cocycle $\phi:\Om\to\Bbb Z^d$ is {\tt D-F} iff $\exists\ N\in\Bbb N,\ M>0$,\ 
$\a_N:=\bigvee_{k=0}^{N-1}\tau^{-k}\a$-measurable functions $g_A:A\to\Bbb Z^d$ ($A\in\a$) with $\|g_A\|_\infty\le M$ and $H:\a\to\Bbb Z^d$ so that
$$\phi|_A=H(A)+g_A\ \forall\ A\in\a.$$
\

For $(\Om,\mu,\tau,\a)$ an {\tt AFU} map, a cocycle $\phi:\Om\to\Bbb Z^d$ is {\tt D-F} iff $\exists\ \ M>0$,\ 
step functions $g_A:A\to\Bbb Z^d$ ($A\in\a$) with $\|g_A\|_\infty\le M$ and $H:\a\to\Bbb Z^d$ so that
$$\phi|_A=H(A)+g_A\ \forall\ A\in\a.$$

\subsection*{Reduction of  cocycles}
\

It is standard that for an {\tt EPPT} $(\Om,\mu,\tau),\ \phi:\Om\to\Bbb G$ measurable, $\chi\in\widehat{\Bbb G},\  \l\in\Bbb S^1\ \& \ g\in L^1(\mu)$:
$$P_{\phi,\chi}(g)=\l g\ \iff\ \chi(\phi)g=\l g\circ\tau$$ and in this case,  $|g|$ is constant. 
\

If in addition, $(\Om,\mu,\tau,\a)$ is a {\tt D-F} fibered system and $\phi:\Om\to\Bbb G$ is  admissible, this situation is 
characterized by $\rho(P_{\phi,\chi})=1$ and entails $g\in L_\tau$.

\

Thus for $\phi:\Om\to\Bbb G$   admissible over a {\tt D-F} fibered system:
\bul $\phi$ is non-arithmetic (as on p. \pageref{Periodicity}) iff for   $\chi\in\widehat{\Bbb G}\setminus\{\mathbb{1}\},\  1$ is not an  eigenvalue of $P_{\phi,\chi}$, and
\bul $\phi$ is aperiodic  iff
for   $\chi\in\widehat{\Bbb G}\setminus\{\mathbb{1}\},\  P_{\phi,\chi}$ has no eigenvalues on the unit circle.

As in \cite[\S3]{AD}, let
$$\frak Q_\phi:=\{\chi\in\widehat{\Bbb G}:\ \rho(P_{\phi,\chi})=1\}=\{\chi\in\widehat{\Bbb G}:\ \exists\ \l\in\Bbb S^1,\ g\in L^1,\ \chi(\phi)g=\l g\circ\tau\},$$ then
\bul \cite[Proposition 3.8]{AD}:  $\frak Q_\phi$ is a closed subgroup of   $\widehat{\Bbb G}$;
\bul $\phi$ is aperiodic iff $\frak Q_\phi=\{\mathbb{1}\}$

Let
$$\Bbb K_\phi=\frak Q_\phi^\perp:=\{x\in\Bbb G:\ \chi(x)=1\ \forall\ \chi\in\frak Q_\phi\}$$
a closed subgroup of $\Bbb G$.

\proclaim{Cocycle reduction formula}\label{reduction}
\

Let $(\Om,\mu,\tau,\a)$ be a {\tt DF} fibered system, either {\tt G-M} or {\tt AFU}.
\

Suppose that $\phi:\Om\to\Bbb G=\Bbb R$ or $\Bbb Z$ is a non-arithmetic, {\tt D-F} cocycle, then
\begin{align*}\tag{\phone}\label{phone}
 \phi=\xi-g+g\circ\tau+\psi
\end{align*}

 where $\xi\in\Bbb G\ \&\ \overline{\Bbb Z\xi+\Bbb K_\phi}=\Bbb G$, $g\in L_\tau\ \&\ \psi:\Om\to\Bbb K_\phi$ is locally $L_\tau$  and aperiodic.
\endproclaim
We'll call the formula  (\phone) the {\it reduction} of $\phi$. 

\

In case $(\Om,\mu,\tau,\a)$ is {\tt G-M}, (\phone)  follows from the topological Markov property of $\tau$ and is a version of  Livsic's cohomology theorem (\cite{Livsic}). See also \cite[Lemma 4.3]{ANS}.
\

\demo{Proof}\ If $\phi$ is aperiodic, then $\frak Q_\phi=\{\mathbb{1}\},\ \Bbb K_\phi=\Bbb G$ and we may set $\xi=0,\ g\equiv 0\ \&\ \psi\equiv \phi$.
\

Otherwise, by non-arithmeticity, $\frak Q_\phi^\perp\cong\Bbb K_\phi\cong\Bbb Z$. 
\

By \cite[Proposition 3.8]{AD}, $\exists\ \xi\in\widehat{\frak Q}_\phi\subset\Bbb G$ and
  $G_\chi\in L_\tau$ ($\chi\in\frak Q$) satisfying $P_{\phi,\chi}(G_\chi)=\xi(\chi)G_\chi$.
  \

  Possibly renormalizing the $G_\chi$, we can ensure that, in addition,  $|G_\chi|=1\ \&\ G_{\chi+\chi'}=G_\chi G_{\chi'}$. 
\

Thus $\exists\ g:\Om\to\Bbb G$ so that $G_\chi=\chi(g)$. Standard lifting theory now shows that $g\in L_\tau$.
\

Lastly, define $\psi:\Om\to\Bbb G$ by
$$\psi:=\phi-\xi+g-g\circ\tau.$$
Evidently $\psi$ is locally $L_\tau$. 
\

Moreover, $\chi(\psi)\equiv 1\ \forall\ \chi\in\frak Q_\phi$, whence $\psi:\Om\to\Bbb K_\phi$.
\

To see that $\overline{\Bbb Z\xi+\Bbb K_\phi}=\Bbb G$, suppose otherwise, then $\exists\ \chi\in\widehat{\Bbb G}\setminus\{\mathbb{1}\}$ so that
  $\chi(n\xi+k)=1\ \forall\ n\in\Bbb Z,\ k\in\Bbb K_\phi$. It follows that
  $$\chi(\phi)\chi(g)=\chi(\xi+g\circ\tau+\psi)=\chi(g\circ\tau)$$
  contradicting non-arithmeticity.
  \

To see that $\psi$ is aperiodic suppose that $\chi\in\widehat{\Bbb G},\ \l\in\Bbb S^1$ and $H\in L^1(\mu)$ satisfy
$\chi(\psi)H=\l H\circ\tau,$ then
\begin{align*}\chi(\phi)H\chi(g)&=\chi(\psi+\xi+g\circ\tau)\\ &=
\chi(\psi)Hg\chi(\xi)\tfrac{\chi(g\circ\tau)}{\chi(g)}\\ &=
\l \chi(\xi)H\circ\tau\chi(g\circ\tau)\end{align*}
whence $\chi\in\frak Q_\phi$ and $\l\chi(\xi)=\chi(\xi)$ entailing $\l=1$.\ \  \Checkedbox

 If $\phi$ is  aperiodic, the aperiodic {\tt LLT}s  follow from \cite[\S6]{AD} ({\tt G-M} case)
 $\&$ \cite[\S5]{ADSZ} ({\tt AFU} case).
 \
 
Our proof of Theorem B is based on this and our advertised
\

\proclaim{Theorem 3.1  {\rm( Periodic {\tt LLT})}}\label{GenLLT} 
\

Suppose that  $(\Om,\mu,\tau,\a)$ is a  {\tt D-F} fibered system, {\tt G-M} or {\tt AFU}. 
\

Let $\phi:\Om\to\Bbb G=\Bbb R$ or $\Bbb Z$  be a positive {\tt D-F} cocycle  so that $(\Om,\mu,\tau,\phi)$ is an  asymptotically $\g$-stable, positive, ergodic, stationary, stochastic process with $0<\g<1$.
\

If $\phi:\Om\to\Bbb G$ is non-arithmetic but not aperiodic,
then $(\Om,\mu,\tau,\a,\phi)$ satisfies 
 a  periodic {\tt LLT}.  \endproclaim

 The theorem  for   independent, identically distributed random variables is  \cite[Theorem 3]{Shepp}. 
 \

 
\

The rest of this section is devoted to the 
\subsection*{Proof of Theorem 3.1}
\

  By  Nagaev's theorem (\cite{Nagaev}), (see also theorem 4.1 in \cite{AD}):
 \

{\it There are constants $\e>0,\ K>0$ and $\th\in (0,1)$; and  continuous functions
$\l=\l_\phi:B(0,\e)\to B_{\Bbb C}(0,1),\  N=N_\phi:B(0,\e)\to \hom(L_\tau,L_\tau)$ such that
$$\|P_{\phi,t}^nh-\l(t)^n N(t)h\|_{L_\tau}\le K\th^n\|h\|_{L_\tau}
 \ \ \ \ \forall\ |t|<\e,\ n\ge 1,\ h\in L_\tau.$$}
Here $P_{\phi,t}:=P_{\phi,\chi_t}$ with $\chi_t(s)=e(ts)$ is as defined on p. \pageref{char} and for $|t|<\epsilon$, $N(t)$ is a projection onto a one-dimensional subspace (spanned by $g(t):=N(t)\mathbb{1}$).

In view of  (\ref{faKey}) on p. \pageref{faKey},  \cite[Theorem 5.1]{AD} applies in the {\tt G-M} case and
\cite[Theorem 8]{ADSZ} applies in the {\tt AFU} case to show that
$$|\l_\phi(t)-\int_\Om \chi_t(\phi)d\mu|=o\left(\frac1{a(\frac1{|t|})}\right)\ \text{as}\ t\to 0.$$
Consequently, 
\begin{align*}\tag{\Biohazard}\label{Biohazard}\l_\phi(\tfrac{t}{a^{-1}(n)})^n\xrightarrow[n\to\infty]{}\ \Bbb E(\chi_t(Z_\g)).
\end{align*}

\demo{Proof of theorem 3.1 in the continuous case}
\

Possibly renormalizing, we assume WLOG that $\Bbb K_\phi=\Bbb Z$ in the cocycle reduction formula
on p. \pageref{reduction} and $\phi$ has reduction 
$$\phi=\xi-F+F\circ\tau+\psi$$
 where $\xi\in (0,1)\ \&\ \overline{\Bbb Z\xi+\Bbb Z}=\Bbb R$, $F\in L_\tau\ \&\ \psi:\Om\to\Bbb Z$ is locally $L_\tau$  and aperiodic.
\

It suffices to show, for fixed $h\in L_\tau$,
 $g\in L^1(\Bbb R)$ with Fourier transform $\widehat{g}\in C^\infty_C(\Bbb R)$ and 
\ \ $k_n\in\Bbb R_+,\ \frac{k_n}{a^{-1}(n)}\xrightarrow[n\to\infty]{}\ \kappa\in\Bbb R_+$,  that
 \begin{align*}\tag{\Bleech}a^{-1}(n)\widehat{\tau}^n(h g(\phi_n-k_n))\approx\ G_g(n\xi-k_n)\Bbb E(h)f_{Z_\g}(\kappa)  
 \end{align*}
 where $G_g(y):=\sum_{J\in\Bbb Z}g(y+J)\ \&\ f_{Z_\g}$ is the probability density of  $Z_\g$.

 \

 Writing $e(t):=e^{2\pi it}$, we have
 \begin{align*}\ttau^n&(hg(\phi_n-k_n))=\int_\Bbb R\widehat{g}(z)e(k_nz)\ttau^n(e(-z\phi_n)h)dz\\ &=
 \int_\Bbb R\widehat{g}(z)e(k_nz)\ttau^n(e(-z(\psi_n+n\xi+F-F\circ\tau^n))h)dz  \\ &=
 \sum_{J\in\Bbb Z}\int_{(-\frac12,\frac12]}\widehat{g}(z+J)e((k_n-n\xi)(z+J))\ttau^n(e(-(z+J)(\psi_n+F-F\circ\tau^n))h)dz. 
 \end{align*}
For $J\in\Bbb Z,\ z\in (-\tfrac12,\tfrac12)$,
 \begin{align*}\widehat{g}(z+J)&e((k_n-n\xi)(z+J))\ttau^n(e(-(z+J)(\psi_n+F-F\circ\tau^n))h)\\ &=
   \widehat{g}(z+J)e((k_n-n\xi)J)e((k_n-n\xi)z)e(-zF)P_{\psi,-z}^n(e(-zF)h)
 \end{align*}
Set $\G_{n,g}(z):= \sum_{J\in\Bbb Z}\widehat{g}(z+J)e((k_n-n\xi)J)$, then by the Poisson summation formula, 
$\G_{g,n}(0)=G_g(k_n-n\xi)$, and

 \begin{align*}\ttau^n&(hg(\phi_n-k_n))=
 \int_{(-\frac12,\frac12]}\G_{n,g}(z)e((k_n-n\xi)z)e(-zF)P_{\psi,-z}^n(e(-zF)h)dz\\ &=
 \int_{-\d }^\d e(z(-F+n\xi-k_n)) \G_{n,g}(z)\l_{\psi}(z)^n N_\psi(z)(he(zF))dz+O(\rho^n) \ \text{in}\ L_\tau
 \\ &=
 \Bbb E(h)\int_{-\d }^\d e(z(n\xi-k_n)) \G_{n,g}(z)\l_{\psi}(z)^n dz+\mathcal E_n+O(\rho^n) 
 \end{align*}

where
$$\mathcal E_n:=\int_{-\d }^\d \G_{n,g}(z)e^{-izk_n}\l_{\phi}(z)^n(e(-zF)N_\psi(z)(he(zF))-\Bbb E(h))dz.$$
We claim that
\begin{align*}\tag{\dsheraldical}\label{dsheraldical}\|\mathcal E_n\|_{L_\tau}=O(\tfrac1{a^{-1}(n)^2})\ \text{as}\ n\to\infty. 
\end{align*}
            
To see this, set $\frak n(z):=e(-zF)N_\psi(z)(he(zF))$, then  by   \cite[theorem 2.4]{AD} ({\tt G-M} case) and \cite[Proposition 5]{ADSZ} ({\tt AFU} case), for some constant $K>0$
$$\|\frak n(z)-\Bbb E(h)\|_{L_\tau}\le K\|h\|_{L_\tau}|z|.$$
Thus
\begin{align*}\|\mathcal E_n\|_{L_\tau}  &\le
\|g\|_1\int_{-\d }^\d |\l_{\phi}(z)|^n\|\frak n(z)-\mathbb{E}(h)\|_{L_\tau}dz\\ &\le K\|g\|_1\|h\|_{L_\tau}\int_{-\d }^\d |\l_{\phi}(z)|^n|z|dz
\end{align*}
 and
\begin{align*}\int_{-\d }^\d |\l_{\phi}(z)|^n|z|dz&=\frac1{a^{-1}(n)^2}\int_{-\d a^{-1}(n)}^{\d a^{-1}(n)} |x||\l_{\phi}(\tfrac{x}{a^{-1}(n)})|^ndx
\\ &\sim \frac1{a^{-1}(n)^2}\int_\Bbb R |x|e^{-C|x|^\g}dx\ \text{with}\ C:=\frac{(2\pi)^\g\cos(\frac{\g\pi}2)}{\G(1+\g)}\\ &
\ \ \ \ \text{by (\Biohazard) on p. \pageref{Biohazard}.\ \Checkedbox\ (\dsheraldical)}
\end{align*}

Thus, 

 \begin{align*}a^{-1}(n)&\widehat{\tau}^n(hg(\phi_n-k_n))\\ &=a^{-1}(n)\Bbb E(h)\int_{-\d }^\d e(z(n\xi-k_n)) \G_{n,g}(z)\l_{\psi}(z)^n dz+O(\tfrac1{a^{-1}(n)})\ \text{in $L_\tau$},\\ &=
 \Bbb E(h)\int_{-\d a^{-1}(n)}^{\d a^{-1}(n)} e(\tfrac{z(n\xi-k_n)}{a^{-1}(n)}) \G_{n,g}(\tfrac{z}{a^{-1}(n)})\l_{\psi}(\tfrac{z}{a^{-1}(n)})^n dz\\ &\approx
 \G_{n,g}(0)\Bbb E(h)\int_\Bbb R e(-z\kappa)\widehat{f}_{Z_\g}(z)dz\\ &=
 G_g(k_n-n\xi)\Bbb E(h)f_{Z_\g}(\kappa)\ \text{by the aperiodic {\tt LLT}.\ \ \ \Checkedbox\ (\Bleech)}
 \end{align*}
 \demo{Proof of theorem 3.1 in the lattice case}
 \

 Suppose that $\tfrac{k_n}{a^{-1}(n)}\xrightarrow[n\to\infty]{}\ \kappa$. We'll show that
   \begin{align*}\tag{\faRocket}a^{-1}(n)\widehat{\tau}^n(h1_{[\phi_n=k_n]})\ \approx\ {p}1_{{p}\Bbb Z}(n\xi-k_n)\Bbb E(h)f_{Z_\g}(\kappa)\ \text{as}\ n\to\infty.    
   \end{align*}

  \begin{align*}\widehat{\tau}^n(&h1_{[\phi_n=k_n]})(x)=\int_\Bbb T e(-k_nt)\ttau^n(e(t\phi_n)h)(x)dt\\ &=
  \int_\Bbb T e(-k_nt)\ttau^n(e(t(p\psi_n+n\xi+F-F\circ\tau^n))h)(x)dt\\ &=
    \int_\Bbb T e(t(n\xi-k_n))\ttau^n(e(p\psi_n)e(t(F-F\circ\tau^n))h)(x)dt\\ & \ \ \ \ \ \ \ \ \ \ 
    \because\ \sum_{q=0}^{p-1}e(\tfrac{Kq}p)=1_{p\Bbb Z}(K)\ \&\ F\circ\tau^n=F\ \mod p;\\ &=
        \int_\Bbb T (\sum_{q=0}^{p-1}(e(\tfrac{t+q}{p}(n\xi-k_n))\ttau^n(e(t\psi_n)e(\tfrac{t+q}{p}(F-F\circ\tau^n))h)(x))dt\\ &=
        {p}1_{{p}\Bbb Z}(n\xi-k_n)\int_\Bbb T e(\tfrac{t}{p}(n\xi-k_n))\ttau^n(e(t\psi_n)e(\tfrac{t}{p}(F-F\circ\tau^n))h)(x))dt\\ &=
      {p}1_{{p}\Bbb Z}(n\xi-k_n)\int_\Bbb T e(\tfrac{t}{p}(n\xi-k_n-F(x))P_{\psi,t}^n(e(\tfrac{tF}{p})h)(x)dt.
   \end{align*}
   As in the proof of the continuous case, using the aperiodic {\tt LLT},
   $$a^{-1}(n)\int_\Bbb T e(\tfrac{t}{p}(n\xi-k_n-F(x))P_{\psi,t}^n(e(\tfrac{tF}{p})h)(x)dt\xrightarrow[n\to\infty]{}\ \Bbb E(h)f_{Z_\g}(\kappa).\ \CheckedBox
   \text{\ (\faRocket)}$$
\section*{\S4 Proof of theorem B}
\subsection*{Uniform sets for $T$ } 

Let the  {\tt MPT} $(X,m,T)$ be pointwise dual ergodic.
\

The set $\Om\in\B(X),\ \ 0<m(\Om)<\infty$  is called 
\bul a {\it uniform set} for $T$ if for some $f\in L^1(m)_+$, the convergence in ({\tt PDE}) (on p. \pageref{PDE}) is uniform on $\Om$; and
\bul  a {\it Darling-Kac set} if this is the case for $f=1_\Om$.
By  Egorov's theorem, every pointwise dual ergodic {\tt MPT} $(X,m,T)$  has uniform sets, which form a  dense hereditary collection   denoted by $\mathcal U(T)$. 

\

If  $(X,m,T)$ is pointwise dual ergodic $\g$-regularly varying return sequence $a(n)$ ($0<\g<1$)., then as in  
\cite[ch. 3]{IET}, the return time process to a uniform set is asymptotically $\g$-stable (defined on p. \pageref{assp}).

\subsection*{ Proof of Theorem B}\ \  Write $\mu=m_\Om,\ \v=\v_\Om\ \&\ \tau=T_\Om$ and let $\a\subset\B(\Om)$ be the partition for which
$(\Om,\mu,\tau,\a)$ is mixing {\tt G-M}/{\tt AFU} map. 
\

If $(\Om,\mu,\tau,\a)$ is {\tt G-M}, then as in   \cite[\S1]{AD}, 
$(\Om,\mu,\tau,\v)$ is continued fraction mixing   and, by the Main Lemma in \cite{RFEx}, $\Om$ is a Darling Kac set for $T$. 
By the asymptotic renewal equation 
(see   e.g. \cite[chapter 3]{IET}))
we  have (\ref{faKey})  (as  on p. \pageref{faKey}).
\

It is shown in \cite{coeffts} that a mixing {\tt AFU} map  $(\Om,\mu,\tau,\a)$ is exponentially reverse $\phi$-mixing.
Thus, by \cite[theorem 2.2]{A-Z}, the pointwise  dual ergodicity of  $(\Om,\mu,\tau)$ with $\g$-regularly varying return sequence ensures (\ref{faKey}).  
\

By our remarks on periodicity of cocycles on p. \pageref{Periodicity}, $\v:\Om\to\Bbb Z$ is non-arithmetic. The required {\tt GLLT} follows from Theorem 3.1 on p.
 \pageref{GenLLT}.
\Checkedbox
\subsection*{Example:  Random walk skew products}\label{RWSP}
\

\ Let $(\Om,\mu,T,\a)$ be a mixing {\tt PP} Gibbs-Markov map and let 
$\phi:\Om\to\Bbb Z$
be $\a$-measurable and aperiodic with $\mu-\text{\tt dist}(\phi)\in\text{\tt DA}\,(SqS)$ with $1<q\le 2$, then
 $(\Om\x\Bbb G,\mu\x\#,T_\phi)$ is:
conservative, exact and pointwise dual ergodic with $a(n)\ \g=1-\tfrac1q$-regularly varying and  $\Om\x\{0\}$ is a uniform set for $T_\phi$. See \cite{AD}.
\

\proclaim{Proposition D}\label{propnD}\ \  $\Om\x\{0\}$ is a  {\tt GLLT} set for $T_\phi$.
\endproclaim\demo{Proof}\  We'll apply Theorem B. To this end, since $(\Om\x\Bbb G,\mu\x\#,T_\phi)$ is exact, hence weakly mixing, it suffices to
exhibit a suitable  Markov partition for the return time process to ${\Om\x\{0\}}$.
\

Evidently  $\v_{\Om\x\{0\}}:\Om\x\{0\}\to\Bbb N$ is given by 
\begin{align*}
 \v_{\Om\x\{0\}}(x,0)&:=\v(x):=\min\,\{n\ge 1:\ \phi_n(x)=0\}\ \&\\ & T_{\Om\x\{0\}}(x,0)=(\tau(x),0):=(T^{\v(x)}(x),0).
\end{align*}

\


Evidently $[\v=n]$ is $\a_n$-measurable. Write $[\v=n]=\bigcup_{a\in\frak b_n}a$ where
$\frak b_n\subset\a_n$ and let $\b\subset\B (\Om)$ be the  partition defined by  
$$\b:=\bigcup_{k\ge 1}\frak b_k.$$
\

It follows that  $(\Om,\mu,\tau,\b)$ is a mixing Gibbs Markov map 
and
 $\v:\Om\to\Bbb N$ is  $\b$-measurable. 
 
The Proposition now follows from  Theorem B.\ \ \Checkedbox
\

\section*{\S5 Continuous time}\label{semiflow}

A measure preserving {\it semiflow} is a continuous semigroup homomorphism $\Psi:\Bbb R_+\to\text{\tt MPT}\,\xbm$  where $\xbm$ is a $\s$-finite, polish measure space and $\text{\tt MPT}\,\xbm$ is the topological semigroup of measure preserving transformations equipped with the weak operator topology.
\

Let 
$(\Om,\mu,\tau)$ be a probability preserving transformation 
and let $\frak r:\Om\to\Bbb R_+$ be measurable.
\

Define the {\it suspended semiflow} $(\Om,\mu,\tau)^\frak r=(X,m,\Psi)$  by
$$X:=\{(x,s)\in\Om\x\Bbb R_+:\ 0\le s<\frak r(x)\},\ m:=\mu\x\text{\tt Leb}$$ and define 
$\Psi:\Bbb R_+\to\text{\tt MPT}\xbm$  by 
$$\Psi_t(x,y)=(\tau^n(x),y+t-\frak r_n(x))$$ where
$n=n_t(x,y)$ is so that
$$0\le y+t-\frak r_n(x)<\frak r(\tau^nx)\ \ \text{\tt\small i.e.}\  y+t\in [\frak r_n(x),\frak r_{n+1}(x)).$$

In this case, we call $(\Om,\mu,\tau,\frak r)$ a {\it section} of the  semiflow $(X,m,\Psi)$, 
which in turn is aka the {\it suspension} of $(\Om,\mu,\tau)$ (under the ceiling   $\frak r$).
\

We'll call the section $(\Om,\mu,\tau,\frak r)$ {\it standard} if the ceiling is bounded below in the sense that $\frak r\ge\D$ a.s. for some $\D=:\min\frak r>0$.
\subsection*{{\tt D-K} sections for  pointwise dual ergodic   semiflows}\label{standard} \ 
Call the semiflow $\Psi$  {\it pointwise dual ergodic}   if each {\tt MPT} $\Psi_t$ ($t>0$) is   pointwise dual ergodic. In this case, 
$$a_t(n):=a_n(\Psi_t)\sim \frac{a(nt)}t$$ (with $a(n):=a_1(n)$.

We'll call a standard section $(\Om,\mu,\tau,\frak r)$ a {\it Darling Kac} ({\tt D-K}) {\it section} for $(X,m,\Psi)=(\Om,\mu,\tau)^\frak r$  if 
$\Om\x [0,t]$ is a Darling Kac set for $\Psi_t\ \forall\ 0<t<\min\frak r$.

\proclaim{Proposition 5.1}\ \ Suppose that $(\Om,\mu,\tau,\frak r)$ is a  {\tt D-K} section for $\Psi$ and
that $a_n(\Psi_1)=:a(n)$ is  $\g$-r.v. ($\g\in (0,1)$), then
\begin{align*}\tag{\Football}\label{Football}
 \mu([\frak r>t])\ \sim\ \frac1{\G(1+\g)\G(1-\g)}\cdot\frac1{a(t)}.
\end{align*}
\endproclaim

\demo{Proof}
\

Fix $0<\D<\min\frak r$, then $\Om_\D:=\Om\x [0,\D]$ is a Darling Kac set for $\Psi_\D$. As shown above, the return sequence 
$a_\D(n):=a_n(\Psi_\D)\sim \tfrac{a(n\D)}\D$.
\

Let $\v=\v_{\Om_\D}:\Om_\D\to\Bbb N$ denote the  first return time to $\Om_\D$ under iterates of $\Psi_\D$,  then for $(x,y)\in\Om_\D=\Om\x [0,\D],\ \v(x,y)=\lcl \tfrac{\frak r(x)-y}\D\rcl$ 
\

By the asymptotic renewal equation 
(see   e.g. \cite[chapter 3]{IET}))
$$\int_{\Om_\D}(\v\wedge t)dm\sim\frac1{\G(1+\g)\G(2-\g)}\cdot\frac{t}{a_\D(t)}$$
and as $t\to\infty$, whence by Karamata's differentiation theorem,
\begin{align*}
 m([\v>t])\sim \frac{d}{dt}\int_\Om(\v_\Om\wedge t)dm\ \sim\ \frac{C_\g}{a_\D(t)}
\end{align*}
with  $C_\g:=\frac1{\G(1-\g)\G(1+\g)}$.

Using this and the sandwich
$$[\v>\tfrac{t}\D]\supset [\frak r>t]\x [0,\D]\supset \{(x,y)\in\Om_\D:\ \frak r(x)-y>t\}\supset [\v>\tfrac{t}\D+1],$$
 we have, 
\begin{align*} \frac{C_\g}{ a(t) }\sim\frac{C_\g}{\D a_\D(\frac{t}\D) }\sim\frac1\D m([\v>\frac{t}\D])\sim\mu([\frak r>t]).\ \ \CheckedBox
\end{align*}

\subsection*{Local limit sections for  pointwise dual ergodic   semiflows}\ \ 

 A a {\it generalized local limit { \tt(GLL)}  section} for $(X,m,\Psi)$ is a quintuple   $(\Om,\mu,\tau,\a,\frak r)$ where
$(X,m,\Psi)=(\Om,\mu,\tau)^{\frak r}$, $(\Om,\mu,\tau,\a)$ is a fibered system as on p. \pageref{fibered} and   $(\Om,\mu,\tau,\a,\frak r)$ satisfies a 
 {\tt GLLT} as defined on p. \pageref{LLT}.

\proclaim{Theorem 5.2}\label{theorem5.2}
\

Suppose that  $(X,m,\Psi)$  is pointwise dual ergodic with $\g$-regularly varying  return sequence {\rm($\g\in (0,1)$)} and which has a  {\tt D-K}, {\tt GLL} section, then
for $t>0$, 
\begin{align*}\tag{\dstechnical}\label{dstechnical}\tfrac1{a(N)}\sum_{n=1}^N&|\widehat{\Psi_t}^n(1_Cg(\tfrac{S_n^{(\Psi_t)}(f)}{a_t(n)}))-m(C) \Bbb E(g(m(f)W_\g))u_t(n)|
\xrightarrow[N\to\infty]{}\ 0\\ & \text{a.e.}\ \forall\ C\in\mathcal F_+,\ g\in C_B(\Bbb R_+)\ \&\ f\in L^1_+ 
\end{align*}
where  $a_t(n):=a_n(\Psi_t),\ u_t(n)\sim\tfrac{\g a_t(n)}n,\ W_\g\in\text{\tt RV}\,(\Bbb R_+),\  \Bbb E(g(W_\g))=\Bbb E(Y_\g g(Y_\g))$.
\endproclaim
The 
property  (\dstechnical)\ \ ($\forall\ t>0$) is a semiflow analogue of (\ref{dsmilitary}) on p. \pageref{dsmilitary}\ \ and is required for \ref{faBatteryQuarter}\ \ (on p. \pageref{faBatteryQuarter}). See Remark 5.4 (on p. \pageref{dsmathematical}) for a ''time-free`` version of \ref{dstechnical}.

\

The proof of Theorem 5.2 mirrors that of   theorem A; beginning with the continuous time version of lemma 2.1: 

\

 \proclaim{Lemma 5.3}
\

Suppose that  $(X,m,\Psi)=(\Om,\mu,\tau)^\frak r$  is pointwise dual ergodic with $\g$-regularly varying  return sequence
$a_n(\Psi_1)=:a(n)$ {\rm($\g\in (0,1)$)} with $(\Om,\mu,\tau,\a,\frak r)$  a periodic, {\tt LLT}, {\tt D-K} section. 
\

Suppose also that $a$ satisfies \ref{faPlug} on p. \pageref{faPlug}, then
for $0<\D<\min\,\frak r,\ A\x I\subset X$ with $A\in\mathcal C_\a\ \&\  I\subset \Bbb R_+$ an interval,

\begin{align*}&\tag{GL*}\label{GL*}\varliminf_{n\to\infty}\tfrac1{u_\D(n)}\widehat{\Psi_{\D}}^n(1_{A\x I}
g(\tfrac{S_n^{(\Psi_\D)}(1_{\Om\x [0,\D]})}{a_\D(n)}))\ge \mu(A)|I| \Bbb E(g(\D W_\g)).
\end{align*} 
\endproclaim
 \demo{Proof of {\rm \ref{GL*}}}\label{lemma5.3}

   Fix $A\in\mathcal C_\a(\tau),\ I\subset \Bbb R_+\ \&\ g\in C_B(\Bbb R_+)$ as above  and fix $0<c<d<\infty$. Assume WLOG that
   $I\subset [pk,p(k+1))$ for some $k\in\Bbb Z$ where $p$ is the period of the {\tt LLT} of the section (p. \pageref{LLT}).
   \
   
   For $N\ge 1$,
   \begin{align*}\widehat{\Psi}^N_{\D}(1_{A\x I}g(\tfrac{S^{(\Psi_\D)}_{N}(1_{\Om\x [0,\D]})}{a_\D(N)}))(x,y) &=   \sum_{n\ge 1}\widehat{\tau}^n(1_{A\cap [\frak r_n\in N\D-y+I]}g(\tfrac{S^{(\Psi_\D)}_{N}(1_{\Om\x [0,\D]})}{a_\D(N)}))\\&=\sum_{n\ge 1}\widehat{\tau}^n(1_{A\cap [\frak r_n\in N\D-y+I]}g(\tfrac{n}{a_\D(N)})).
\end{align*}
   
   Writing $x_{n,N}:=\tfrac{N\D}{a^{-1}(n)}$, we have by 
 regular variation that
 \begin{align*}   \tag{{\scriptsize\faLeaf}}\label{faLeaf}
   n\sim \tfrac{\D^\g a(N)}{x_{n,N}^\g}\ &\&\ 
\tfrac1{a^{-1}(n)}\sim (\tfrac1{x_{n,N}})^\g\cdot (x_{n,N}-x_{n+1,N})\cdot\tfrac{u(N)}{\D^{1-\g}}\\ &  \text{as}\ \ N,\ n\to\infty,\ x_{n,N}\in [c,d].
 \end{align*}

Let $$K_N:=\{n\ge 1,\ x_{n,N}\in [c,d]\}.$$
Using the {\tt LLT} property of $\Om$, (\ref{faLeaf}) and Lemma 2.2, we have, as $N\to\infty$,
\begin{align*}\tag{{\scriptsize\faBug}}\label{faBug}&\widehat{\Psi}^N_{\D}(1_{A\x I}g(\tfrac{S^{(\Psi_\D)}_{N}(1_{\Om\x [0,\D]})}{a_\D(N)}))(x,y)\\ & \ge \sum_{n\in K_N}g(\tfrac{n}{a_\D(N)})\widehat{\tau}^n(1_{A\cap [\frak r_n\in N\D-y+I]})\\ &\sim
\sum_{n\in K_N}g(\tfrac{\D}{x_{n,N}^\g})\widehat{\tau}^n(1_{A\cap [\frak r_n\in N\D-y+I]})\\ &=
\sum_{n\in K_N}g(\tfrac{\D}{x_{n,N}^\g})\widehat{\tau}^n(1_{A\cap [\frak r_n\in x_{n,N}a^{-1}(n)-y+I]})\\ &\sim
\sum_{n\in K_N}g(\tfrac{\D}{x_{n,N}^\g})\tfrac{pf_{Z_\g}(x_{n,N})\mu(A)}{a^{-1}(n)}
1_{p\Bbb Z+I}(N\D-n\xi)                                  
\\ &\sim                                                                  
\tfrac{u(N)\mu(A)}{\D^{1-\g}}\sum_{n\in K_N}g(\tfrac{\D}{x_{n,N}^\g})\tfrac{p f_{Z_\g}(x_{n,N})}{x_{n,N}^\g}(x_{n,N}-x_{n+1,N})
1_{p\Bbb Z+I}(N\D-n\xi)\\ & 
\sim \tfrac{u(N)\mu(A)}{\D^{1-\g}}\Bbb E(1_{[c,d]}(Z_\g)g(\D Z_\g^{-\g})Z_\g^{-\g})\ \cdot\ pm_{\Bbb R/p\Bbb Z}(I+p\Bbb Z)\\ &=
u_\D(N)\mu(A)|I| \Bbb E(1_{[c,d]}(Z_\g)g(\D Z_\g^{-\g})Z_\g^{-\g}).\ \ \ \CheckedBox\ \text{\ref{GL*}}
\end{align*}

\demo{Completion of the proof  Theorem 5.2 (p. \pageref{theorem5.2})}
\ 

The proof now proceeds as the completion of the proof of theorem A on p. \pageref{A-finish} to obtain \ref{dstechnical}for 
$0<t<\min\frak r$. This suffices as if $\Psi_t$ satisfies (\dstechnical), then so does $\Psi_{Nt}=\Psi_t^N\ \forall\ N\ge 1$.\ \Checkedbox
\

\subsection*{Remark 5.4:  A ``{\tt\small time-free}'' consequence of \ref{dstechnical}}
\

\begin{align*}\tag{\dsmathematical}\label{dsmathematical}\tfrac1{a(N)}\int_1^N
&|\widehat{\Psi_t}(1_Cg(\tfrac{\Xi_t^{(\Psi)}(f)}{a(t)}))-m(C) \Bbb E(g(m(f)W_\g))u(t)|dt\\ &
\xrightarrow[N\to\infty]{}\ 0\ \text{a.e.}\ \forall\ C\in\mathcal F_+,\ g\in C_B(\Bbb R_+)\ \&\ f\in L^1_+ 
 \end{align*}
 where $\Xi_t^{(\Psi)}(f):=\int_0^tf\circ\Psi_sds$.
 \demo{Proof sketch} The idea here is to use the identity:
 $$\Xi_n^{(\Psi)}(f):=\int_0^nf\circ\Psi_sds=S_n^{(\Psi_1)}(\overline{f})\ \text{with}\ \overline{f}:=\int_0^1f\circ\Psi_sds.$$ 
We have, for $f\in L^1_+\cap L^\infty,\ \log g\in C([0,\infty])_+,\ n\ge 1\ \&\ r\in (0,1)$
  \begin{align*}\Xi_{n+r}^{(\Psi)}(f)&=\Xi_{n}^{(\Psi)}(f)\circ\Psi_r-\int_0^rf\circ\Psi_sds\\ &=
  S_{n}^{(\Psi_1)}(\overline{f})\circ\Psi_r\pm \overline{f}.
     \end{align*}
Thus, uniformly as $n\to\infty$, 
$$g(\tfrac{\Xi_{n+r}^{(\Psi)}(f)}{a(n+r)})\sim
g(\tfrac{S_{n}^{(\Psi_1)}(\overline{f})}{a(n)})\circ\Psi_r$$
and for $C\in L^1_+\cap L^\infty$,
 \begin{align*}\tag*{{\scriptsize\faPlane}}\label{faPlane}
 \int_n^{n+1}
\widehat{\Psi_t}(Cg(\tfrac{\Xi_t^{(\Psi)}(f)}{a(t)}))dt&=
\int_0^1\widehat{\Psi}_{n+r}(Cg(\tfrac{\Xi_{n+r}^{(\Psi)}(f)}{a(n+r)}))dr\\ &\sim
\int_0^1\widehat{\Psi}_{n+r}(Cg(\tfrac{S_{n}^{(\Psi_{1})}(\overline{f})\circ\Psi_r}{a(n)}))dr\\ &=
\widehat{\Psi}_{n}(\widetilde{C}\cdot g(\tfrac{S_{n}^{(\Psi_{1})}(\overline{f})}{a(n)}))\end{align*}

with $\widetilde{C}:=\int_0^1\widehat{\Psi}_{r}(C)dr.$
 
To use this, extend   (with a similar proof) \ref{dstechnical} with $t=1$ to:
 \begin{align*}\tag{\dsagricultural}\label{dsagricultural}\tfrac1{a(N)}\sum_{n=1}^N&|\widehat{\Psi_1}^n(Cg(\tfrac{S_n^{(\Psi_1)}(f)}{a(n)}))-m(C) \Bbb E(g(m(f)W_\g))u(n)|
\xrightarrow[N\to\infty]{}\ 0\\ & \text{a.e.}\ \forall\ C\in L^1_+,\ g\in C_B(\Bbb R_+)\ \&\ f\in L^1_+
\end{align*}
where $m(C):=\int_XCdm$ for $C\in L^1_+$.
\

By \ref{faPlane} and \ref{dsagricultural} for a.e. $(x,y)\in X$, there is a set $K_{(x,y)}\subset\Bbb N$ of full density so that
\begin{align*}
\tfrac1{u(n)}\int_n^{n+1}
\widehat{\Psi_t}(&Cg(\tfrac{\Xi_t^{(\Psi)}(f)}{a(t)}))(x,y)dt\approx \tfrac1{u(n)}\widehat{\Psi}_{n}(\widetilde{C}\cdot g(\tfrac{S_{n}^{(\Psi_{1})}(\overline{f})}{a(n)}))(x,y)\\ & \xrightarrow[n\to\infty,\ n\in K_{(x,y)}]{}\ m(C) \Bbb E(g(m(f)W_\g)),
 \end{align*}
whence \ref{dsmathematical}.\ \CheckedBox

\subsection*{Remark 5.5}\label{remark5.5}
In certain cases, it is possible to treat non-standard sections (see p. \pageref{standard})  in theorem 5.2 by obtaining a standard {\tt\small induced section}.
\

Let $(X,m,T)=(\Om,\mu,\tau)^\frak r$ be pointwise dual ergodic with $\g$-regularly varying return sequence $a(n)$ ($\g\in (0,1)$, and suppose that $\frak r\ge\D>0$ on $A\in\B(\Om)$. It is well known that
$(X,\mu(A)^{-1}m,T)$ is a factor semiflow of $(A,\mu_A,\tau_A)^{\widehat{\frak r}}$ where
$$\widehat{\frak r}(x):=\sum_{k=0}^{\v_A(x)-1}\frak r(\tau^kx)$$
with $(A,\mu_A,T_A,\v_A)$  the  return time process to $A$. 
\

Moreover $(A,\mu_A,\tau_A,\widehat{\frak r})$ is a standard section of $(A,\mu_A,\tau_A)^{\widehat{\frak r}}$.
\

If $(\Om,\mu,\tau,\frak r)$ is a positive, Bernoulli process, then so is $(A,\mu_A,\tau_A,\widehat{\frak r})$. 
The latter is also standard, whence a {\tt D-K}, {\tt GLLT} section.
\

Similarly, if $(\Om,\mu,\tau,\a)$ is a {\tt G-M} map, $\frak r:\Om\to\Bbb R_+$ is locally $L_\tau$ and $A\in\a$, then
$(\Om,\mu_A,\tau_A,\b)$ is also a {\tt G-M} map with
$$\b:=\bigcup_{n=1}^\infty\{a\in\a_n:\ \v|_a\equiv n\}$$
and $\widehat{\frak r}$ is locally $L_{\tau_A}$. Again, it follows that $(A,\mu_A,\tau_A,\widehat{\frak r})$ is
 a {\tt D-K}, {\tt GLLT} section.

\subsection*{ Remarks about mixing}\label{sfHKmixing} \ \ \
\

 Let      $\Psi=(\Om,\mu,\tau)^\frak r$ be a pointwise, dual ergodic semiflow with $a(n)=a_n(\Psi_1)$ $\g$-regularly varying ($\g\in (0,1)$). Suppose that
 $(\Om,\mu,\tau,\a,\frak r)$ is a {\tt GLLT} section.
 \
 
    In a similar manner to the remarks about mixing on p. \pageref{HKmixing}, it follows from \ref{GL*} $\&$ (\ref{faBug}) on p. \pageref{faBug} that   TFAE:
\begin{align*}&\tag{i}\tfrac1{u_\D(N)}\widehat{\Psi}^N_{\D}(1_{A\x I}g(\tfrac{S^{(\Psi_\D)}_{N}(1_{\Om\x [0,\D]})}{a_\D(N)}))\xrightarrow[N\to\infty]{} \mu(A)|I| \Bbb E(g(\D W_\g)) 
\\ &  \ \ \text{\small a.e. on $\Om$}
\ \ \forall\ A\in\mathcal C_\a(\tau),\ I\subset (0,\D)\ \text{\small an interval},\ g\in C_B(\Bbb R),\ g\ge 0;
\\ &\tag{ii}\varlimsup_{N\to\infty}\tfrac1{u_\D(N)}
\sum_{n\ge 1,\ x_{n,N}\notin [c,d]}\widehat{\tau}^n(1_{A\cap [\frak r_n\in x_{n,N}a^{-1}(n)-y+[0,\D ]})\xrightarrow[c\to 0+,\ d\to\infty]{}\ 0;\\ &\tag{iii}
\tfrac1{u_\D(N)}\widehat{\Psi}^N_{\D}(1_{\Om\x [0,\D]})\xrightarrow[t\to\infty]{} \D.
\end{align*}

\

\subsection*{Random walk semiflows}\label{RWsf}
\

Let $(\Om,\mu,T,\a)$ be a Gibbs-Markov map as defined on \pageref{GibbsMarkov}, let  $\frak h:\Om\to\Bbb R_+$ be $(\a,\th)$-H\"older as defined on p. \pageref{Holder}, and let $\phi:\Om\to\Bbb Z$ be with $\a$-measurable and aperiodic with $\mu-\text{\tt dist}(\phi)\in\text{\tt DA}\,(SqS)$ with $1<q\le 2$.

\

Consider the measure preserving semiflow 
$$(X,m,\Psi)=(\Om\x\Bbb Z,\mu\x\#,T_\phi)^{\overline{\frak h}}$$  where:
\bul $(\Om\x\Bbb Z,\mu\x\#,T_\phi)$ is the corresponding 
random walk skew product as on p. \pageref{RWSP}; and  
\bul $\overline{\frak h}(x,z)=\frak h(x)$.

\

Recall that $(\Om\x\Bbb Z,\mu\x\#,T_\phi)$ is 
 exact, pointwise dual ergodic with $a(n)=a_n(T)$ $\g$-regularly varying with $\g=1-\tfrac1q$,  
conditional {\tt RWM} with rate $u(n)=\tfrac{\g a(n)}n$ and $\Om\x\{0\}$ is a {\tt GLL} set for $T_\phi$.

\

We show that, under certain conditions, $\Psi$ is  pointwise dual ergodic  and 
has a  {\tt GLL} section, whence has the tied-down $\g$-renewal mixing property \ref{dstechnical}   (on p. \pageref{dstechnical}).

\

Suppose that $(\frak h,\phi):\Om\to\Bbb R\x\Bbb Z$ is aperiodic, then:
\bul $\Psi$ is conditionally {\tt RWM} with rate $\propto u(t)$ {\small\rm(Theorem 2  in \cite{A-T})};
\bul $\Om\x [0,\D)$ is a uniform set for $\Psi_\D$ for $0<\D<\min\,\frak h$.
\

We conclude by showing that $\Psi$ satisfies \ref{dstechnical}in   theorem 5.2. This is done by finding a  {\tt GLL} section for $\Psi$.

As in the proof of Proposition D on p. \pageref{propnD}, $(\Om,\mu,\tau,\b)$ is a mixing Gibbs Markov map where
\bul $\v:\Om\to\Bbb N$ by $\v(x):=\min\,\{n\ge 1:\ \phi_n(x)=0\}$;
\bul $\tau:\Om\to\Om$ by $\tau(x):=T^{\v(x)}(x)$;
\bul $\b\subset\B (\Om)$ a partition by  $\b=\{a\in\a_n:\ \v|_a\equiv n\}$
and  the induced map of $(T_\phi)_{\Om\x\{0\}}\cong\tau$: 
$$(T_\phi)_{\Om\x\{0\}}(x,0)=(\tau(x),0).$$ 
As above
 $(\Om,\mu,\tau,\b)$ is a mixing Gibbs Markov map.
\

The required section for $\Psi$ will be
$$(\Om,\mu,\tau,\b,\frak r)$$
where
$$\frak r(x):=\sum_{k=0}^{\v(x)-1}\frak h(T^kx).$$
Calculation shows that $D_{\b,\th}(\frak r)<\infty$ and the weak mixing of $\Psi$ ensures that $\frak r$ is non-arithmetic.
It follows as above that $(\Om,\mu,\tau,\b,\frak r)$ is a  {\tt GLL} section for $\Psi$,
and the claim follows from  theorem 5.2.\ \Checkedbox 
\section*{\S6 Integrated properties}

\subsection*{Weak rational ergodicity}\ \
\

The {\tt CEMPT} $(X,m,T)$ is called {\it weakly rationally ergodic} ({\tt WRE}) (\cite{RE}) 
if $\exists$   $F\in\mathcal F_+$ so  that
\begin{align*}\tag{$\largestar$}\frac1{a_n(F)}\sum_{k=0}^{n-1}m(B\cap T^{-k}C)\underset{n\to\infty}\lra\ m(B)m(C)\ \forall\ B,C\in\B\cap F\end{align*}
where $a_n(F):=\frac1{m(F)^2}\sum_{k=0}^{n-1}m(F\cap T^{-k}F)$. 
\

By theorem 3.3 in \cite{FL}, $F\in\mathcal F_+$ satisfies ($\largestar$) if and only if
\begin{align*}\{S_n^{(T)}(1_F):\ n\in\Bbb N\}\ \ \text{\tt\small is uniformly integrable on}\ F.\end{align*}
A useful sufficient condition for this (\cite{RE}, \S3.3 in \cite{IET}) is
$$\sup_{n\ge 1}\tfrac1{a_n(F)^2}\int_FS_n(1_F)^2dm<\infty$$  and  $(X,m,T)$ is called {\it rationally ergodic}  if $\exists$  such  $F\in\mathcal F_+$. 
\

In case $T$ is weakly rationally ergodic:
\sbul the collection of sets $R(T)$ satisfying ($\largestar$)  is  a hereditary ring;
\sbul  $\exists\ a_n(T)$ (the {\it return sequence}) $\st$
$$\frac{a_n(A)}{a_n(T)}\underset{n\to\infty}\lra 1\ \forall\ A\in R(T);$$
\sbul for conservative, ergodic $T$,  $R(T)=\mathcal F$ only when $m(X)<\infty$.
\

By Proposition 3.7.1 in \cite{IET} a pointwise dual ergodic transformation is rationally ergodic with the same 
return sequence.
\subsection*{Rational weak mixing}\ \

As in \cite{RatWM}, we call the {\tt CEMPT} $(X,m,T)$   {\it rationally weakly mixing} ({\tt RWM}) if 
$\exists$ rates $u_n>0$ and $F\in\mathcal F_+$ so that
\begin{align*}\tag{$\bigstar$}\frac1{a(n)}\sum_{k=0}^{n-1}|&m(A\cap T^{-k}B)-u_k(F)m(A)m(B)|
 \\ &\xrightarrow[n\to\infty]{}\ 0
\ \ \forall\ A,\ B\in \B\cap F\end{align*}
where $a(n):=\sum_{k=0}^{n-1}u_k$.
\

It is shown in \cite{RatWM} that {\tt RWM} entails {\tt WRE} with ($\bigstar$) holding $\forall\ F\in R(T)$.
\proclaim{Proposition 6.1}\ If $(X,m,T)$ is conditionally {\tt RWM}, then it is {\tt RWM}.\endproclaim
Proposition 6.1 is implicit in \cite{nice}. The method of proof is the same as that of the next result.

 \proclaim{Theorem 6.2}\ \ Suppose that  $(X,m,T)$ satisfies {\rm (\dsmilitary)} on p. \pageref{dsmilitary}, then
$(X,m,T)$ is {\tt WRE} and
\begin{align*}\tag{\scriptsize\faTrain}\label{faTrain}\tfrac1{a(N)}&\sum_{n=1}^N|\int_{B\cap T^{-n}C}g(\tfrac{S_n(f)}{a(n)})dm-m(B)m(C)
 \Bbb E(g(m(f)W_\g))u(n)|\\ &
\xrightarrow[N\to\infty]{}\ 0\ \forall\ B,C\in\mathcal R\cap R(T),\ g\in C_B(\Bbb R_+)\ \&\ f\in L^1_+
\end{align*} where $a(n)=a_n(T)$.\endproclaim
The integrated, tied-down renewal mixing property   {\rm \scriptsize(\faTrain)} of the {\tt MPT} $(X,m,T)$ is a strengthening of  {\tt RWM}  (take $g\equiv 1$ in  {\rm \scriptsize(\faTrain)}).

\

The proof of Theorem 6.2 uses a standard approximation technique embodied in
\proclaim{Lemma 6.3} Suppose that $(X,m,T)$ is {\tt WRE},  $f\in L^1_+,\ g\in C_B(\Bbb R_+)_+,\ \Om\in R(T)$  and that {\rm \scriptsize(\faTrain)}
holds for $A\in\mathcal R,\ B\in\mathcal S$ where both $\mathcal R,\ \mathcal S\subset\B(\Om)$ are dense in $\B(\Om)$, then
{\rm \scriptsize(\faTrain)}  
holds $\forall\ A, B\in\B(\Om)$.\endproclaim
\demo{Proof}\ 

We claim first that  {\rm \scriptsize(\faTrain)} holds for $A\in\B(\Om)\ \&\ b\in\mathcal S$. 
\

Indeed for $A\in\B(\Om)\ \&\ \e>0,\ \exists\ a\in\mathcal R$ so that $m(a\D A)<\e$ whence
\begin{align*}|&\int_{A\cap T^{-k}b}g(\tfrac{S_n(f)}{a(n)})dm-u_km(A)m(bu)\Bbb E(m(f)W_\g)|\\ &\le 
|\int_{a\cap T^{-k}b}g(\tfrac{S_n(f)}{a(n)})dm-u_km(a)m(bu)\Bbb E(m(f)W_\g)|+\\ &\ \ \ +\|g\|_{C_B}(m(a\D A\cap T^{-k}b)+u_km(a\D A)m(b)).  
 \end{align*}

Using {\rm \scriptsize(\faTrain)} for $a\in\mathcal R,\ b\in\mathcal S$ and $\Om\in R(T)$, we have as $n\to\infty$
 \begin{align*}
\tfrac1{a(n)}\sum_{k=0}^{n-1}&|\int_{A\cap T^{-k}b}g(\tfrac{S_n(f)}{a(n)})dm-u_km(A)m(bu)\Bbb E(m(f)W_\g)|\lesssim\\ &\tfrac{\|g\|_{B_B}}{a(n)}\sum_{k=0}^{n-1}(m(a\D A\cap T^{-k}b)+u_km(a\D A)m(b))\\ &\sim 2 m(a\D A)m(b)<2\e.
 \end{align*}\
 The extension of  {\rm \scriptsize(\faTrain)} to $A,\ B\in\B(\Om)$ is similar.\ \Checkedbox
\demo{Proof of Theorem 6.2}\ Fix $\Om\in \mathcal R_0:=R(T)\cap\mathcal R$, then $a_n(\Om)\sim a(n)$.  
Fix and $f\in L^1(m)_+$.
\

Using Egorov's theorem, for fixed $A\in \B(\Om)$ and $\e>0$, $\exists\ U=U_A\in\B(\Om)$ so that $m(\Om\setminus U)<\e$ and so that
$$\tfrac1{a(n)}\sum_{k=1}^n|\T^k(1_Ag(\tfrac{S_n(f)}{a(n)}))-u_km(A)\Bbb E(m(f)W_\g)|\xrightarrow[n\to\infty]{}0\ \text{uniformly on}\ U,$$
whence for $B\in\B(U)$,
\begin{align*}
 \frac1{a(n)}\sum_{k=0}^{n-1}&|\int_{A\cap T^{-k}B}g(\tfrac{S_n(f)}{a(n)})dm-u_km(A)m(B)\Bbb E(m(f)W_\g)|\\ &=
 \frac1{a(n)}\sum_{k=0}^{n-1}|\int_B(\T^k(1_Ag(\tfrac{S_n(f)}{a(n)}))-u_km(A)\Bbb E(m(f)W_\g))dm|\\ &\le
 \int_B(\tfrac1{a(n)}\sum_{k=0}^{n-1}|\T^k(1_Ag(\tfrac{S_n(f)}{a(n)}))-u_km(A)\Bbb E(m(f)W_\g)|)dm\\ &\xrightarrow[n\to\infty]{}0.
\end{align*}
Now fix a countable dense algebra $\mathcal R\subset\B(\Om)$. 
\

By the above, $\exists\ \{U_{A,n}:\ n\ge 1,\ A\in\mathcal R\}\subset\B(\Om)$ so that
 $m(V_n)>m(\Om)-\frac1{n}$ where $V_n:=\bigcap_{A\in\mathcal R}U_{A,n}$ and so that
 \begin{align*}
\tag*{{\rm \scriptsize(\faTrain)}}\frac1{a(n)}\sum_{k=0}^{n-1}&|\int_{a\cap T^{-k}b}g(\tfrac{S_n(f)}{a(n)})dm-u_km(a)m(bu)\Bbb E(m(f)W_\g)|\xrightarrow[n\to\infty]{}0\\ & \forall\ a\in\mathcal R\ \&\ 
b\in\mathcal S:=\bigcup_{n\ge 1}\B(V_n).
 \end{align*}
The extension of  {\rm \scriptsize(\faTrain)} to $A,B\in\B(\Om)$ follows from Lemma 6.2. 
\ \Checkedbox

\subsection*{Natural extensions}
\par Suppose that
$(X,m,T)$ is a measure preserving transformation of a standard,
$\sigma$-finite measure space. Rokhlin's {\it natural extension} $(X^\Bbb N,\widetilde{m},\widetilde{T})$
the minimal, invertible extension of $(X,m,T)$ (which is unique up to isomorphism) is given (as in \cite{NatExt}) by
\begin{align*}
& \widetilde{T}(x_1,x_2,...):=(Tx_1,x_1,...),\\ &
 \widetilde{m}([A_1,...,A_n])=m(\bigcap_{k=1}^nT^{-(n-k)}A_k)\  \ \text{where}\ \\ &
 [A_1,...,A_n]=\{(x_1,x_2,..)\in X^\Bbb N:x_k\in A_k\ \forall 1\le k\le n\},\ 
A_1,...,A_n\in\B(X).
\end{align*}
This is a {\tt MPT} extending $(X,m,T)$ via $\pi:(x_1,x_2,..)\mapsto x_1$. 
\

It is invertible on 
$\widetilde{X}=\{(x_1,x_2,..)\in X^\Bbb N:Tx_{n+1}=x_n\ \forall n\in\Bbb N\}$ which has full measure; with inverse $\widetilde{T}^{-1}(x_1,x_2,..)=(x_2,x_3,\dots)$. Thus
$(\widetilde{X},\widetilde{m},\widetilde{T})$ is an invertible {\tt MPT}.
\

The extension is minimal since
$\bigvee_{n=0}^\infty \widetilde{T}{}^n\pi^{-1}\B(X)=\B(\widetilde{X})\ \mod\ \widetilde{m}.$

We say that a property $\frak P$ (of {\tt MPT}s) {\it lifts } (to the natural extension) if
$$(X,m,T)\ \text{has property}\ \frak P\ \Rightarrow\ (\widetilde{X},\widetilde{m},\widetilde{T}) \ \text{has property}\ \frak P.$$

Parry showed in \cite{ErgSpec} that conservativity and ergodicity  lift.
It follows (see \cite{IET}) that rational ergodicity and weak rational ergodicity both lift and it is shown in
\cite{RatWM} that rational weak mixing lifts. The following shows that the integrated tied-down renewal mixing property (\ref{faTrain}) on p. \pageref{faTrain} also lifts.
\proclaim{Theorem 6.4}\ \ \ \ Suppose that$(X,m,T)$ is {\tt WRE} with $a(n)=a_n(T)$  $\g$-regularly varying with  $\g\in (0,1)$ and satisfies {\rm(\ref{faTrain})}, then so does $(\widetilde{X},\widetilde{m},\widetilde{T})$.
\endproclaim\demo{Proof} Let $\pi:(\widetilde{X},\widetilde{m},\widetilde{T})\to (X,m,T)$ be the extension map and write $a(n)=a_n(T)\sim a_n(\widetilde{T})$.

\

Let $\mathcal R$ be the hereditary ring in the $\g$-tied-down renewal mixing of  $(X,m,T)$ and set 
$$\mathcal R^*:=\bigcup_{n\in\Bbb Z}\widetilde{T}^{n}\pi^{-1}\mathcal R.$$
We have that $\mathcal R^*\subset R(\widetilde{T})$ and that 
for $\Om\in\mathcal R^*,\ \mathcal R^*\cap\B(\Om)$ is dense in $\B(\Om)$.
\

Also, for fixed $f\in L^1(m)_+,\ g\in C_B(\Bbb R_+)_+$,
{\small\begin{align*}\tag*{{\rm\scriptsize(\faBus)}}\tfrac1{a(N)}&\sum_{n=1}^N|\int_{B\cap \widetilde{T}^{-n}C}g(\tfrac{S_n(f\circ\pi)}{a(n)})d\widetilde{m}-\widetilde{m}(B)\widetilde{m}(C) \Bbb E(g(\widetilde{m}(f)W_\g))u(n)|\\ &
\xrightarrow[N\to\infty]{}\ 0\ \forall\ B,C\in\mathcal R^*.
\end{align*}}
Thus by Lemma 6.2, for each $\Om\in\mathcal R^*$, ({\rm\scriptsize\faBus}) holds 
$\forall\ B,C\in\B(\Om)$ and, putting things together, we see that ({\rm\scriptsize\faBus}) holds
\sms   $\forall\ f=p\circ\pi,\ p\in L^1(m)_+,\ g\in C_B(\Bbb R_+)\ \&$
 $B,C\in\widetilde{\mathcal R}:=\bigcup_{\Om\in\mathcal R^*}\B(\Om)$. 

\

We now extend  this to general $F\in L^1(\widetilde{m})_+$. The method is similar to the proof of theorem A.
\

Again, it suffices to consider $g>0$ with
 $x\mapsto\ \log g(e^x)$  uniformly continuous on $\Bbb R$. 
 
Fix $F,\ f\in L^1_+(\widetilde{m})_+$, $f=p\circ\pi,\ p\in L^1(m)_+,\ m(p)=1,\ g>0$ as above and $\Om\in\mathcal R^*$.

\

By the ratio ergodic theorem
$$\frac{S_n(F)}{S_n(f)}\xrightarrow[n\to\infty]{}\ m(F) \ \text{a.e.}$$
By Egorov's theorem
$$\mathcal U:=\{C\in\B(\Om):\ \text{convergence uniform on}\ C\}$$ is dense in $\B(\Om)$.
\

Let $B\in\mathcal U$, then $\exists\ \e_N\downarrow 0$ so that
$$g(\tfrac{S_n(F)}{a(n)})=(1\pm\e_n)g(\tfrac{m(F)S_n(f)}{a(n)})\ \text{on}\ B\ \forall\ n\ge 1.$$
Thus for $C\in\B(\Om)$,
\begin{align*}|\int_{B\cap\widetilde{T}^{-n}C} g(&\tfrac{S_n(F)}{a(n)})d\widetilde{m}-\widetilde{m}(B)\widetilde{m}(C) \Bbb E(g(\widetilde{m}(F)W_\g))u(n)|\\ & \le 
|\int_{B\cap\widetilde{T}^{-n}C} (g(\tfrac{S_n(F)}{a(n)})-g(\tfrac{\widetilde{m}(F)S_n(f)}{a(n)}))|d\widetilde{m}+\\ &\ \ \ \
|\int_{B\cap\widetilde{T}^{-n}C} g(\tfrac{\widetilde{m}(F)S_n(f)}{a(n)})d\widetilde{m}-\widetilde{m}(C) \Bbb E(g(\widetilde{m}(F)W_\g))u(n)|\\ &\le
\e_n\|g\|_\infty\widetilde{m}(B\cap\widetilde{T}^{-n}C)+\\ &
\ \ \ +|\int_{B\cap\widetilde{T}^{-n}C} g(\tfrac{\widetilde{m}(F)S_n(f)}{a(n)})d\widetilde{m}-\widetilde{m}(C) \Bbb E(g(\widetilde{m}(F)W_\g))u(n)|.\end{align*}
Now
$$\tfrac1{a(N)}\sum_{n=1}^N|\int_{B\cap\widetilde{T}^{-n}C} g(\tfrac{\widetilde{m}(F)S_n(f)}{a(n)})-\widetilde{m}(C) \Bbb E(g(\widetilde{m}(F)W_\g))u(n)|\xrightarrow[N\to\infty]{}0,$$
so
\begin{align*}\tfrac1{a(N)}\sum_{n=1}^N&|\int_{B\cap\widetilde{T}^{-n}C} g(\tfrac{S_n(F)}{a(n)})-\widetilde{m}(B)\widetilde{m}(C) \Bbb E(g(\widetilde{m}(F)W_\g))u(n)|\\ &\lesssim\ 
\tfrac1{a(N)}\sum_{n=1}^N\e_n\|g\|_\infty\widetilde{m}(B\cap\widetilde{T}^{-n}C)\\ &
\xrightarrow[N\to\infty]{}\ 0.
\end{align*}
Thus ({\rm\scriptsize\faBus}) holds $\forall\ B\in\mathcal U,\ C\in\B(\Om)$. Lemma 6.2 now extends ({\rm\scriptsize\faBus})  to hold $\forall\ B,\ C\in\B(\Om)$.\ \Checkedbox
\subsection*{Example: Geodesic flows on cyclic covers }\label{cyclic}
\

We denote by $(U(M),\L,g)$, the {\tt\small geodesic flow} on the {\tt\small unit tangent bundle} $U(M)$ of the  {\tt\small hyperbolic surface} $M$ equipped with
the hyperbolic measure $\L$ (for definitions see \cite{Hopf},\ \cite{Rees},\ \cite{nice}). 
\

The hyperbolic surface $V$ is a {\it cyclic cover} of the compact  hyperbolic surface $M$ if there is
a covering map $p:V\to M$ and a monomorphism $\g:\Bbb Z\to\text{\tt Isom}\,(V)$ (hyperbolic isomotries of $V$), so that
for $y\in V,\ p^{-1}\{p(y)\}=\{\g(n)y:\ n\in\Bbb Z\}$. 
\

Symbolic dynamics for the geodesic flow on a compact hyperbolic surface is described in \cite{Bowen} and it is shown in \cite{Rees} (see also \cite{nice}) that if $V$ is a  cyclic cover of a compact  hyperbolic surface, then
$(U(M),\L,g)$ is isomorphic to the natural extension of a semiflow of form $(X\x\Bbb Z,m\x \#,T_\varphi)^{\overline{\frak r}}$
where $(X,T,m,\a)$ is a mixing Gibbs Markov map with $\#\a<\infty$ (also known as a subshift of finite type equipped with a Gibbs measure) and $(\phi,\frak r):X\to\Bbb Z\x\Bbb R$ is $\a$-H\"older and (\cite{Solomyak})
aperiodic. It was shown in \cite{nice} that such $(U(M),\L,g)$ is rationally weakly mixing. For more information, see
\cite{A-T}.
\

We claim here that such 
\Smi $(U(M),\L,g)$ satisfies (\ref{faTrain}) on p. \pageref{faTrain} with $a(n)\propto\sqrt n$.

\demo{Proof of \smiley} \ \ As shown in the random walk semiflow example  on p. \pageref{RWsf}, 
$(X\x\Bbb Z,m\x \#,T_\varphi)^{\overline{\frak r}}$  (\ref{faTrain})  with $a(n)\propto\sqrt n$.
\smiley\  now follows from theorems 6.2 and 6.4.\ \Checkedbox
\subsection*{Updates}\ \ A functional version of \ref{faTrain} is established in \cite{A-Sera-functional} under stronger assumptions; as is a $u$-{\tt weak Cesaro} (as in \cite[\S4]{RatWM}) version for pointwise dual ergodic transformations  with $\g$-regularly varying
return sequences ($0<\g\le 1$). This latter gives versions of \ref{faLightbulbO} $\&$ \ref{faBatteryQuarter}  (p.
\pageref{faLightbulbO}) when $\g=1$.

\end{document}